\documentclass[11pt,letterpaper]{amsart}
\usepackage{amsmath,amssymb,amsfonts,amsthm}
\usepackage[normalem]{ulem} 
\usepackage{soul} 
\usepackage{mathrsfs}
\usepackage[all]{xy}
\usepackage{stmaryrd}
\usepackage{mathtools}
\usepackage{indentfirst}
\usepackage{verbatim}
\usepackage{hyperref}
\usepackage{color}
\usepackage[usenames,dvipsnames]{xcolor}
\usepackage{graphicx}
\graphicspath{ {images/} }
\usepackage{tikz}
\usepackage{nicefrac}
\usepackage{xfrac}
\usepackage[shortlabels]{enumitem}

\title{Log canonical models of foliated surfaces}
\author{Yen-An Chen}
\thanks{The author was partially supported by NSF research grants no: DMS-1801851, DMS-1840190 and by a grant from the Simons Foundation; Award Number: 256202. }
\address{National Center for Theoretical Sciences, Taipei City 106, Taiwan (R.O.C.)}
\email{yachen@ncts.tw}

\newtheorem*{claim}{Claim}
\newtheorem{thm}{Theorem}[section]
\newtheorem{prop}[thm]{Proposition}
\newtheorem{cor}[thm]{Corollary}
\newtheorem{lem}[thm]{Lemma}

\theoremstyle{definition}
\newtheorem{defn}[thm]{Definition}
\newtheorem*{acks}{Acknowledgements}
\theoremstyle{remark}
\newtheorem{rmk}[thm]{Remark}
\newtheorem*{pf}{Proof}

\newtheorem{eg}[thm]{Example}

\newcommand{\supp}{\mbox{{\rm Supp}}}

\newcommand{\proj}{\mbox{{\rm Proj \!}}}

\newcommand\cA{{\mathcal{A}}}

\newcommand\cD{{\mathcal{D}}}
\newcommand\cE{{\mathcal{E}}}
\newcommand\cF{{\mathcal{F}}}
\newcommand\cG{{\mathcal{G}}}

\newcommand\cL{{\mathcal{L}}}
\newcommand\cM{{\mathcal{M}}}

\newcommand\cO{{\mathcal{O}}}

\newcommand\cQ{{\mathcal{Q}}}

\newcommand\cS{{\mathcal{S}}}

\newcommand\cU{{\mathcal{U}}}

\newcommand\cX{{\mathcal{X}}}
\newcommand\cY{{\mathcal{Y}}}
\newcommand\cZ{{\mathcal{Z}}}

\newcommand\bA{{\mathbb A}}

\newcommand\bC{{\mathbb C}}

\newcommand\bN{{\mathbb N}}
\newcommand\bP{{\mathbb P}}
\newcommand\bQ{{\mathbb Q}}
\newcommand\bR{{\mathbb R}}

\newcommand\bZ{{\mathbb Z}}

\newcommand\sF{{\mathscr{F}}}
\newcommand\sG{{\mathscr{G}}}
\newcommand\sH{{\mathscr{H}}}
\newcommand\sL{{\mathscr{L}}}

\newcommand\ord{\mbox{\rm ord}}
\newcommand\rw{\longrightarrow}
\newcommand\wt[1]{\widetilde{#1}}

\newcommand\vol{\mbox{{\rm vol}}}
\newcommand\ve{\varepsilon}

\newcommand{\p}[1]{\frac{\partial}{\partial #1}}
\newcommand{\D}{\Delta}
\newcommand{\T}{\Theta}
\newcommand\tn[1]{\textnormal{#1}}

\begin{document}

\maketitle

\begin{abstract}
We study log canonical models of foliated surfaces of general type. 
In particular, we show that log canonical models of general type and their minimal partial du Val resolutions are bounded. 
Moreover, we show the valuative criteria of separatedness and properness and a property related to local-closedness for the moduli functor $\cM_P^{\tn{sm}}$ which parametrizes the stable smoothable foliated surfaces pairs. 

On the way, we also show a result on the invariance of plurigenera. 
\end{abstract}

\section*{Introduction}
The study of foliated surfaces of general type has seen significant progress in recent years. 
By work of McQuillan and Brunella (see \cite{mcquillan2008canonical}, \cite{brunella2015birational}, and the reference therein), it is known that foliated surfaces of general type with only canonical foliation singularities admit a unique canonical model. 
Thus it is natural to study the moduli theory for these canonical models, for instance, the existence of the moduli functor. 

A first step toward this was achieved in \cite{hacon2021birational} where it is shown that canonical models of foliated surfaces of general type with fixed Hilbert function are birationally bounded. 
Notice that, by a result of McQuillan, $K_\sF$ is not $\bQ$-Cartier for any canonical model $(X,\sF)$ with cusp singularities. 
Thus, the moduli functor is expected not to be algebraic. 
However, if we consider the moduli functor for some natural birational modifications of these canonical models, then it is better behaved. 
For instance, in \cite{chen2021boundedness}, it is shown that minimal partial du Val resolutions (see Definition~\ref{mpdvr}) of these canonical models are bounded. 
As a consequence, the boundedness of these canonical models also follows. 

The next natural property to investigate is the local-closedness. 
We recall the following example. 
\begin{eg}[{\cite[page 114]{brunella2001invariance}}]\label{key_eg}
Consider $\pi : \bC^2_{x,\,y}\times(\bC_{\lambda}\setminus\{0,1\}) \rw \bC_\lambda\setminus\{0,1\}$ and a family of foliations given by 
\[x\p{x}+\lambda y\p{y}.\]
We can extend this family to a family $\bP^2\times(\bC\setminus\{0,1\})\rw\bC\setminus\{0,1\}$. 
Note that the singularity is reduced (see Definition~\ref{semi_reduced}) if and only if $\lambda\not\in\bQ^+$. 
\end{eg}

In light of this example, the moduli functor for canonical models (or their minimal partial du Val resolutions) is \emph{not} locally closed. 
To solve this issue, we allow some special \emph{log canonical} foliation singularities and consider \emph{log canonical} models (see Definition~\ref{defn_gcm}). 
As in \cite{hacon2021birational} and \cite{chen2021boundedness}, we first show the following boundedness for minimal partial du Val resolutions of log canonical models. (See Theorem~\ref{key_thm}.)
\begin{thm}
Fix a function $P: \bZ_{\geq 0} \rw \bZ$. 
Then the set $\cS_P$ consisting of minimal partial du Val resolutions $(Y,\sG,D)$ (see Definition~\ref{mpdvr}) of log canonical models $(X,\sF)$ of general type with Hilbert function $\chi(mK_\sF) = P(m)$ is bounded. 
More precisely, there exists a projective morphism $\mu : \cY \rw T$ and a family of foliations $\cF$ on $\cY$ over $T$ (see Definition~\ref{family_fol}) where $\cY$ and $T$ are quasi-projective varieties of finite type such that, for any $(Y,\sG,D)\in\cS_P$, there is a $t\in T$ and an isomorphism $\rho : Y\rw \cY_t$ such that $\sG(-D) \cong \rho^*(\cF\vert_{\cY_t})$.
\end{thm}

Using the techniques in \cite{cascini2018invariance}, we also show the following invariance of plurigenera for families of minimal partial du Val resolutions of log canonical models of general type with fixed Hilbert function. (See Theorem~\ref{inv}.)
Note that condition $(\cX,\cF,\Delta)_t\in\cS_P$ puts a very strong restriction on $\Delta$. 
\begin{thm}
Given any family of foliated triples $f : (\cX,\cF,\Delta)\rw T$ where $T$ is connected 
and $(\cX,\cF,\Delta)_t\in\cS_P$ for all $t\in T$ (see Definition~\ref{defn_fol_family} for families of foliated triples), there exists an integer $m_{f}$ depending only on this family $f$ such that  
\[h^0(\cX_t,m(K_{\cF}+\Delta)\vert_{\cX_t})\] is independent of $t$ for $m\geq m_f$ and divisible by $i(P)$. (See Lemma~\ref{iP}.)
\end{thm}

Next, for a fixed function $P : \bZ_{\geq 0} \rw \bZ$, we use the techniques from \cite{spicer2020higher} and \cite{cascini2021mmp}, to establish the semi-stable $(K_\cF+\D)$-minimal model program for the semi-stable families of foliated triples. (See Section 4 and Theorem~\ref{KF_MMP}.)

\begin{thm}
Let $f : (\cX,\cF,\D) \rw T$ be a family of foliated triples which is a semi-stable reduction of a family of foliated triples whose general fiber is in $\cS_P$. 
(See subsection~\ref{ss_reduction}.)
Then the relative $(K_\cF+\D)$-minimal model program exists and stops with a family of foliated triples $(\cX^{\tn{m}},\cF^{\tn{m}},\Delta^{\tn{m}}) \rw T$ such that $(K_{\cF^{\tn{m}}}+\D^{\tn{m}})$ is nef over $T$, $f^{\tn{m}} : \cX^{\tn{m}} \rw T$ is dlt, and the general fiber is in $\cS_P$. 
\end{thm}

We then define the moduli functor $\cM_P^{\tn{sm}}$ (see Definition~\ref{functor}) and show it satisfies the valuative criteria for separatedness and properness by utilizing the semi-stable $(K_\cF+\D)$-minimal model program. (See Theorem~\ref{separatedness} and Theorem~\ref{properness}.) 

\begin{thm}[Separatedness]
Let $T$ be a smooth curve, $T^\circ = T\setminus\{0\}$ where $0\in T$, and $(\cX_i,\cF_i,\D_i) \in\cM_P^{\tn{sm}}(T)$ be two families of stable $\cS_P$-smoothable foliated triples over $T$ for $i=1$, $2$. 
Then any isomorphism $\varphi^\circ : (\cX_1,\cF_1,\D_1)_{T^\circ} \rw (\cX_2,\cF_2,\D_2)_{T^\circ}$ over $T^\circ$ extends uniquely to an isomorphism $\varphi : (\cX_1,\cF_1,\D_1) \rw (\cX_2,\cF_2,\D_2)$ over $T$. 
\end{thm}

\begin{thm}[Properness]
Let $T$ be a smooth curve and $(\cX^\circ,\cF^\circ,\D^\circ)\in\cM_P^{\tn{sm}}(T^\circ)$ be a family of stable $\cS_P$-smoothable foliated triples over $T^\circ$ where $T^\circ = T\setminus\{0\}$ and $0\in T$. 
Suppose the fiber of the family over $t\neq 0$ is in $\cS_P$. 
Then there exists a smooth curve $C$ and a finite morphism $C \rw T$ such that, after the base change via $C \rw T$, the pullback family $(\cX_C^\circ,\cF_C^\circ,\D_C^\circ)$ extends uniquely to a family of stable $\cS_P$-smoothable foliated triples $(\cX_C,\cF_C,\D_C)$ over $C$. 
\end{thm}

Finally, by a Seidenberg-type Theorem \cite[Proposition 7.4]{pereira2019effective} for a family of foliated surfaces, we show the following theorem (see also Theorem~\ref{open}), which gives us a hint of the local-closedness of the moduli functor $\cM_P^{\tn{sm}}$. 

\begin{thm}
Given a family of foliated triples $f : (\cX,\cF,\D)\rw T$. 
($T$ is not necessarily a curve.) 
Assume that there is a dense subset $T'\subset T$ such that $(\cX,\cF,\D)_t\in\cS_P$ for all $t\in T'$. 
(See Definition~\ref{set_MPdVR} and Definition~\ref{defn_fol_family}.) 
Then there exists an open subset $U\subset T$ such that all fibers of $f$ over $U$ are in $\cS_P$. 
\end{thm}

\begin{acks}
A large part of the paper was written during the author was at the University of Utah. He would like to thank the University of Utah for wonderful research conditions. The author would also like to thank Christopher D. Hacon for his insightful suggestions and encouragements and Fabio Bernasconi for invaluable discussions. 
\end{acks}

\section{Preliminaries}
In this paper, we always work over $\bC$. 
Throughout the paper, a scheme will refer to a scheme of finite type over $\bC$ and a variety is a complex analytic space. 

In this section, we recall several definitions and results which will be used later. 

\subsection{Intersection theory on normal complete surfaces}
We will use Mumford's intersection pairing on normal surfaces (see, e.g., \cite{sakai1984weil}) and some of its properties. 

Let $X$ be a normal complete surface and $f: Y\rw X$ be a proper birational morphism from a smooth surface $Y$, for instance, the minimal resolution of singularities. 
Let $E = \cup_jE_j$ be the exceptional divisor of $f$. 
Since the intersection matrix $(E_i\cdot E_j)$ is negative definite, for any Weil $\bR$-divisor $D$ on $X$, we can define $f^*D$ as $\wt{D} + \sum_ja_jE_j$ where $\wt{D}$ is the proper transform of $D$ and the $a_j$'s are real numbers uniquely determined by the property that $(\wt{D} + \sum_ja_jE_j)\cdot E_i = 0$ for every $i$. 

For any two Weil $\bR$-divisors $D_1$ and $D_2$ on $X$, we define the intersection number $D_1\cdot D_2 := (f^*D_1)\cdot (f^*D_2)$. 
Note that if $D_1$ is Cartier and $D_2$ is an irreducible curve, then $D_1\cdot D_2$ agrees with the degree of $\cO_{D_2}(D_1)$. 

\subsection{A criterion for very-ampleness}
\begin{thm}[{\cite[Theorem 0.2]{langer2001adjoint}}]\label{va_eff}
Let $X$ be a normal projective surface and $M$ be a $\bQ$-divisor on $X$ such that $K_X+\lceil M\rceil$ is Cartier. 
Let $\zeta$ be a 0-dimensional subscheme of $X$ and $\delta_\zeta>0$ be the number defined in \cite[Definition 1.1.2]{langer2001adjoint} associated to $\zeta$. If 
\begin{enumerate}
\item $M^2>\delta_\zeta$ and 
\item $M\cdot C\geq \frac{1}{2}\delta_\zeta$ for every curve $C$ on $X$,  
\end{enumerate}
then $H^0(X,\cO_X(K_X+\lceil M\rceil)) \rw \cO_\zeta(K_X+\lceil M\rceil)$ is surjective. 
\end{thm}

\begin{rmk}\label{delta}
For the case when $|\zeta| = 2$, we have following remarks. 
\begin{enumerate}
\item If the support of $\zeta = \{x_1,x_2\}$ is two distinct points, then $\delta_{\zeta} = \delta_{x_1}+\delta_{x_2}$. (\cite[1.1.3]{langer2001adjoint})
\item $\delta_x \leq 4$ if $x$ is a smooth point or a du Val singularity. (\cite[0.3.2]{langer2001adjoint})
\item If $\zeta$ is supported at only one point $x$, then 
\begin{enumerate}
\item $\delta_\zeta \leq 8$ if $x$ is a smooth point or a du Val singularity. (\cite[0.3.2]{langer2001adjoint})
\item $\delta_\zeta \leq 4(\textnormal{edim}_xX-1)$ if $x$ is a rational (non-smooth) singularity where $\textnormal{edim}_xX$ is the embedding dimension of $X$ at $x$. (\cite[0.3.4]{langer2001adjoint})
\end{enumerate}
\end{enumerate}
\end{rmk}

\begin{prop}[{\cite[Proposition II.7.3]{hartshorne1977algebraic}}]\label{sep}
Let $X$ be a projective scheme over $k$ and $\sL$ be a line bundle on $X$. 
Then $\sL$ is very ample if and only if it separates points and tangents.
\end{prop}

\subsection{Foliations of rank one on schemes}
\begin{defn}\label{defn_fol_sch}
Let $X$ be a scheme. 
A \emph{pre-foliation} $\sF$ of rank one on $X$ is a rank-one quotient sheaf $\Omega_{\sF}$ of $\Omega_X^{[**]}$. 
Here $\Omega_X^{[**]}$ is the hull of $\Omega_X$ defined in \cite[Definition 14]{kollar2008hulls}. 
A pre-foliation $\sF$ is called a \emph{foliation} if it is \emph{pure}, that is, the closure of each associated point of $\Omega_\sF$ has dimension $\dim X$. 
\end{defn}

\begin{lem}\label{normal_fol_equiv}
Let $X$ be a normal variety. 
A foliation $\sF$ of rank one on $X$ is equivalent to a saturated rank one coherent subsheaf $T_\sF$ of the tangent sheaf $T_X$. 
Here a subsheaf $T_\sF$ of $T_X$ is called saturated if $T_X/T_\sF$ is torsion-free. 
\end{lem}
\begin{pf}
Since $X$ is normal, we have $\Omega_X^{[**]} = \Omega_X^{**}$ the double dual of $\Omega_X$. 

Given a foliation $\sF$ on $X$, we have an exact sequence 
\[\xymatrix{0\ar[r] & K \ar[r] & \Omega_X^{**} \ar[r]^-{\rho_\sF} & \Omega_\sF \ar[r] & 0}\]
where $K$ is the kernel of $\rho_\sF$. 
Dualizing this sequence, we get 
\[\xymatrix{0 \ar[r] & \Omega_\sF^* \ar[r]^-{\rho_\sF^*} & T_X \ar[r] & A \ar[r] & 0}\]
where $A$ is the cokernel of $\rho_\sF^*$. 
Note that $A$ is a subsheaf of $K^*$ and $K^*$ is torsion-free. 
So $A$ is also torsion-free, and thus $\Omega_\sF^*$ is saturated. 
Therefore, this gives a function $F$ from the set of foliations of rank one on $X$ to the set of saturated rank one coherent subsheaves of the tangent sheaf $T_X$. 

Conversely, given a saturated coherent subsheaf $T_\sF$ of $T_X$, we have an exact sequence 
\[\xymatrix{0 \ar[r] & T_\sF \ar[r] & T_X \ar[r]^-\pi & T_X/T_\sF \ar[r] & 0.}\] 
Dualizing the sequence, we have 
\[\xymatrix{0\ar[r] & (T_X/T_\sF)^* \ar[r]^-{\pi^*} & \Omega_X^{**} \ar[r] & B \ar[r] & 0}\]
where $B$ is the cokernel of $\pi^*$. 
Notice that $B$ is a subsheaf of $T_\sF^*$ and $T_\sF^*$ is torsion-free. 
Thus $B$ is also torsion-free and therefore gives a foliation on $X$. 
Hence, this gives a function $G$ from the set of saturated rank one coherent subsheaves of the tangent sheaf $T_X$ to the set of foliations of rank one on $X$. 

Notice that both $G\circ F$ and $F\circ G$ are generically identity functions. 
\begin{claim}
They are indeed identity functions up to isomorphism. 
\end{claim}
\begin{pf}
Notice that we have the following exact sequences: 
\[\xymatrix{0 \ar[r] & K \ar[r] & \Omega_X^{**} \ar[r] \ar@{=}[d] & \Omega_\sF \ar[r] & 0 \\
0 \ar[r] & A^* \ar[r] & \Omega_X^{**} \ar[r] & G(F(\Omega_\sF)) \ar[r] & 0.}\]
Since $G(F(\Omega_\sF))$ is torsion-free and the map 
\[\xymatrix{K\ar[r] & \Omega_X^{**} \ar[r] & G(F(\Omega_\sF))}\] 
is generically zero, the image is in fact zero in $G(F(\Omega_\sF))$. 
Thus, we have $K\subset A^*$. 
Similarly, we have $A^* \subset K$. 
Therefore, we have $K = A^*$ and hence $\Omega_\sF \cong G(F(\Omega_\sF))$. 

Similarly, we could show that $T_\sF = F(G(T_\sF))$. 
\qed
\end{pf}
\end{pf}

\begin{defn}[Foliated triple]\label{fol_triple}
A \emph{foliated triple} is a triple $(X,\sF,D)$ where $X$ is a scheme, $\sF$ is a foliation, and $D$ is a reduced divisor. 
When $D=\emptyset$, we call $(X,\sF) := (X,\sF,\emptyset)$ a \emph{foliated scheme} (resp. \emph{foliated variety}) if $X$ is a scheme (resp. variety). 
\end{defn}

\subsection{Foliations on normal varieties}
Let us recall the following common definition for foliations. 
\begin{defn}
A \emph{foliation} $\sF$ on a normal variety $X$ is a coherent subsheaf $\sF$ of the tangent sheaf $T_X$ such that 
\begin{enumerate}
\item it is closed under the Lie bracket and 
\item $\sF$ is saturated, that is, the quotient $T_X/\sF$ is torsion free.
\end{enumerate}
\end{defn}
The \emph{rank} of $\sF$ is its rank as a sheaf. 
We define the \emph{canonical divisor} $K_\sF$ for the foliation $\sF$ as $\cO_X(-K_\sF) = \det(\sF)$. 

\begin{rmk}
By Lemma~\ref{normal_fol_equiv}, this is equivalent to Definition~\ref{defn_fol_sch} when the rank is one. 
\end{rmk}

Let $\sF$ be a foliation of rank $r$ on a normal variety $X$. 
We have a morphism $\Omega_X^{[1]}\rw \sF^*$. 
Taking the double dual of $r$-th wedge product, we get a morphism 
\[\varphi : \Omega_X^{[r]} \rw \cO_X(K_\sF),\]
which yields a map 
\[\phi : \Omega_X^{[r]}\otimes\cO_X(-K_\sF) \rw \cO_X.\]
We define the \emph{singular locus} of $\sF$ as the co-support of the image of $\phi$. 

\begin{defn}[Rational transform of foliations]
Let $\sF$ be a foliation on a normal variety $X$. 
\begin{enumerate}
\item  Given $f : Y \dashrightarrow X$ be a dominant rational map. 
Let $U$ be an open subset of $X$ such that $f\vert_V : V \longrightarrow U$ is an isomorphism where $V:= f^{-1}(U)$. 
Note that $\sF\vert_U \subset T_U \cong T_V$. 
By \cite[Exercise II.5.15]{hartshorne1977algebraic}, there is a coherent subsheaf $\sG$ of $T_Y$ such that $\sG\vert_V = \sF\vert_U \subset T_V $. 
Then the \emph{pullback foliation} $f^*\sF$ is defined to be the saturation of $\sG$. 
Note that $\sG$ is closed under the Lie bracket and, by \cite[Lemma 1.8]{hacon2021birational}, this definition is well-defined. 
\item Given $g : X \dashrightarrow Z$ be a birational map, then we define the \emph{pushforward foliation} $g_*\sF$ as $(g^{-1})^*\sF$. 
\end{enumerate}
\end{defn}

\begin{defn}
Let $\sF$ be a foliation on a normal variety $X$. 
Given a subvariety $W$ of $X$. 
Let $U = X\setminus\left(\tn{Sing}(X)\cup\tn{Sing}(\sF)\cup\tn{Sing}(W)\right)$. 
\begin{enumerate}
\item $W$ is \emph{tangent} to $\sF$ if $T_W\vert_U\rw T_X\vert_U$ factors through $\sF\vert_U$. 
Otherwise we say $W$ is \emph{transverse} to $\sF$. 
\item $W$ is \emph{invariant} if $\sF\vert_U\rw T_X\vert_U$ factors through $T_W\vert_U$. 
\end{enumerate}
\end{defn}

\begin{defn}
We define Kodaira dimension of $\sF$ as 
\[\kappa(\sF) = \max\{\dim\phi_{|mK_\sF|}(X)\vert\,m\in\bN_{>0} \mbox{ and } h^0(X,mK_\sF)>0\}\] 
where by convention $\kappa(\sF)=-\infty$ if the set above is empty, that is, $h^0(X,mK_\sF)=0$ for all $m\in\bN_{>0}$. 
We say that $\sF$ is \emph{of general type} if $\kappa(K_\sF) = \dim X$. 
\end{defn}

\begin{defn}[{\cite[Definition I.1.5]{mcquillan2008canonical}}]
Let $(\sF,\Delta)$ be a foliated pair on a normal variety $X$ and $f : Y\rw X$ be a proper birational morphism. 
For any divisor $E$ on $Y$, we define the \emph{discrepancy} of $(\sF,\Delta)$ along $E$ to be $a(E,\sF,\Delta) = \textnormal{ord}_E(K_{f^*\sF}-f^*(K_\sF+\Delta))$. 

We say $(X,\sF,\Delta)$ is \emph{terminal} (resp. \emph{canonical}) if $a(E,\sF,\Delta)> 0$ (resp. $\geq 0$) for every exceptional divisor $E$ over $X$ and 
$(X,\sF,\Delta)$ is \emph{log terminal} (resp. \emph{log canonical}) if $a(E,\sF,\Delta)> -\ve(E)$ (resp. $\geq -\ve(E)$) for every divisor $E$ over $X$ where $\varepsilon(E)$ is defined to be $0$ if $E$ is $f^*\sF$-invariant, and $1$ otherwise. 
\end{defn}

\begin{defn}\label{index}
Let $(X,\sF)$ be a foliated normal variety. 
We have the following definitions for indices.
\begin{enumerate}
\item The index $i(X)$ of $X$ is the smallest positive integer $m$ such that $mD$ is Cartier for every Weil divisor $D$ on $X$. 
(We set $i(X)=\infty$ if there is no such $m$.)
\item The index $i(K_X)$ of $K_X$ is the smallest positive integer $m$ such that $mK_X$ is Cartier. 
(We set $i(K_X)=\infty$ if $K_X$ is not $\bQ$-Cartier.)
\item The index $i(\sF)$ of a foliation $\sF$ on $X$ is the smallest positive integer $m$ such that $mK_\sF$ is Cartier. 
(We set $i(\sF)=\infty$ if $K_\sF$ is not $\bQ$-Cartier.)
\item The $\bQ$-index $i_\bQ(\sF)$ of a foliation $\sF$ on $X$ is the smallest positive integer $m$ such that $mK_\sF$ is Cartier for all $\bQ$-Gorenstein points of the foliation. 
\end{enumerate}
\end{defn}

\subsection{Invariants on foliated normal surfaces}
In this subsection, we will focus on some indices on foliated normal surfaces. 
Most definitions follow from \cite{brunella2015birational} with some generalizations. 
Let $(X,\sF)$ be a foliated normal surface. 

Given any $p\in\textnormal{Sing}(\sF)\setminus\textnormal{Sing}(X)$. That is, $p$ is a smooth point on $X$ but a singular point of the foliation $\sF$. 
Let $v$ be the vector field around $p$ generating $\sF$. 
Since $p\in\textnormal{Sing}(\sF)$, we have $v(p)=0$. 
Then we can consider the eigenvalues $\lambda_1$, $\lambda_2$ of $(\tn{D}v)\vert_p$, which do not depend on the choices of $v$. 
\begin{defn}\label{semi_reduced}
If one of the eigenvalues is non-zero, say $\lambda_2$, then we define the eigenvalue of the foliation $\sF$ at $p$ to be 
\[\lambda := \frac{\lambda_1}{\lambda_2}.\] 
For $\lambda\neq 0$, this definition is unique up to reciprocal $\lambda \sim \frac{1}{\lambda}$. 

If $\lambda = 0$, then $p$ is called a \emph{saddle-node}; otherwise, we say $p$ is \emph{non-degenerate}. 
If $\lambda\not\in\bQ^+$, then $p$ is called a \emph{reduced} singularity of $\sF$. 
\end{defn}

\begin{thm}[Seidenberg's theorem]
Given any foliated surface $(X,\sF)$ with $X$ smooth. 
Then there is a sequence of blowups $\pi : (Y,\sG) \rw (X,\sF)$ such that $(Y,\sG)$ has only reduced singularities.
\end{thm}
\begin{pf}
See \cite{seidenberg1968reduction}, \cite[Appendix]{mattei1980holonomie}, or \cite[Theorem 1.1]{brunella2015birational}. \qed
\end{pf}

\begin{defn}[Separatrices]
Let $p$ be a singular point of $\sF$. 
A \emph{separatrix} of $\sF$ at $p$ is a holomorphic (possibly singular) irreducible $\sF$-invariant curve $C$ on a \emph{neighborhood} of $p$ which passes through $p$. 
\end{defn}

\begin{thm}\label{separatrix}
Let $\sF$ be a foliation on a surface $X$ and $C\subset X$ be a connected, compact $\sF$-invariant curve such that the followings hold.
\begin{enumerate}
\item All the singularities of $\sF$ on $C$ are reduced and on the non-singular locus of $X$. 
\item The intersection matrix of $C$ is negative definite, and the dual graph is a tree.
\end{enumerate}
Then there exists at least one point $p\in C\cap\textnormal{Sing}(\sF)$ and a separatrix through $p$ not contained in $C$. 
\end{thm}
\begin{pf}
See \cite{camacho1988quadratic}, \cite{sebastiani1997existence}, or \cite[Theorem 3.4]{brunella2015birational}. \qed
\end{pf}

\subsubsection{Non-invariant curves}
We first consider the non-invariant curves and define the tangency order for them. 
\begin{defn}[Tangency order]
Let $(X,\sF)$ be a foliated normal surface and $p\in\textnormal{Sing}(\sF)\setminus\textnormal{Sing}(X)$. 
Let $v$ be the vector field generating $\sF$ around $p$ and $C$ be a \emph{non-invariant} curve through $p$. 
Let $f$ be the local defining function of $C$ at $p$.  
We define the \emph{tangency order} of $\sF$ along $C$ at $p$ to be 
\[\textnormal{tang}(\sF,C,p) := \dim_\bC\frac{\cO_{X,p}}{\langle f,v(f)\rangle}.\]
Note that $\textnormal{tang}(\sF,C,p)\geq 0$ and is independent of the choices of $v$ and $f$. 
Moreover, if $\sF$ is transverse to $C$ at $p$, then $\textnormal{tang}(\sF,C,p) = 0$. 
Therefore, if $C$ is compact, then we can define 
\[\textnormal{tang}(\sF,C) := \sum_{p\in C}\textnormal{tang}(\sF,C,p). \]
\end{defn}

\begin{thm}[{Adjunction for non-invariant divisors, \cite{brunella1997feuilletages}, \cite[Proposition 2.2]{brunella2015birational}}]\label{adjunct_non_inv}
Let $\sF$ be a foliation on a smooth projective surface $X$. 
Let $C$ be a non-invariant irreducible curve on $X$ and $C^\nu$ be the normalization of $C$.  
Then there is an effective divisor $\Delta$ on $C^\nu$ such that $(K_\sF+C)\vert_{C^\nu} = \Delta$ with $\textnormal{deg}\,\Delta = \textnormal{tang}(\sF,C)$.
\end{thm}

\begin{cor}
Let $C$ be a non-invariant curve on a foliated surface $(X,\sF)$. If $C$ is contained in the smooth locus of $X$, then we have $(K_\sF + C)\cdot C \geq 0$.
\end{cor}

\subsubsection{Invariant curves} 
Now we consider the invariant curves. 
\begin{defn}
Let $(X,\sF)$ be a foliated surface and $p\in\textnormal{Sing}(\sF)\setminus\textnormal{Sing}(X)$. 
Let $\omega$ be a $1$-form generating $\sF$ around $p$. 
If $C$ is an invariant curve and $f$ is the local defining function at $p$, then we can write \[g\omega = h\tn{d}f+f\eta\] where $g$ and $h$ are holomorphic functions, $\eta$ is a holomorphic $1$-form, and $h, f$ are relatively prime functions. 

We define the index $\textnormal{Z}(\sF,C,p)$ to be the vanishing order of $\frac{h}{g}\vert_C$ at $p$. 
Also we define Camacho-Sad index $\textnormal{CS}(\sF,C,p)$ to be the residue of $\frac{-1}{h}\eta\vert_C$ at $p$.  
These two definitions are independent of the choices of $f$, $g$, $h$, $\omega$, $\eta$. 
(For a reference, see \cite[page 15 in Chapter 2 and page 27 in Chapter 3]{brunella2015birational}.)
\end{defn}

Note that if $p\not\in\textnormal{Sing}(\sF)$, then $\textnormal{Z}(\sF,C,p) = 0 = \textnormal{CS}(\sF,C,p)$. 
Therefore, if $C$ is compact, then we can define 
\begin{align*}
\textnormal{Z}(\sF,C) &:= \sum_{p\in C} \textnormal{Z}(\sF,C,p) \tn{ and }\\
\textnormal{CS}(\sF,C) &:= \sum_{p\in C} \textnormal{CS}(\sF,C,p). 
\end{align*}

\begin{thm}[{Adjunction for invariant divisors, \cite{brunella1997feuilletages}, \cite[Proposition 2.3]{brunella2015birational}}]\label{adjunct_inv}
Let $\sF$ be a foliation on a smooth projective surface $X$. 
Let $C$ be an invariant irreducible curve on $X$ and $C^\nu$ be its normalization. 
Then there is an effective divisor $\Delta$ on $C^\nu$ such that $K_\sF\vert_{C^\nu} = K_{C^\nu}+\Delta$ with $\deg\Delta = \textnormal{Z}(\sF,C) + \deg\textnormal{Diff}_C(0)$ where $\textnormal{Diff}_C(0)$ is the different with $(K_X+C)\vert_{C^\nu} = K_{C^\nu}+\textnormal{Diff}_C(0)$. 
In particular, we have $K_\sF\cdot C = \textnormal{Z}(\sF,C) + 2p_a(C)-2$ where $p_a(C)$ is the arithmetic genus of $C$. 
\end{thm}

\begin{thm}[Camacho-Sad formula, \cite{camacho1982invariant}, \cite{linsneto1988algebraic}, \cite{suwa1998indices}]\label{thm_CS_formula}
Let $\sF$ be a foliation on a smooth projective surface $X$, and $C$ be an invariant curve on $X$. 
Then we have $C^2 = \textnormal{CS}(\sF,C)$.
\end{thm}

\begin{lem}\label{lem_ZCS}
Given a foliated normal surface $(X,\sF)$. 
Let $p\in\tn{Sing}(\sF)\setminus\tn{Sing}(X)$ and $\omega$ be a $1$-form generating $\sF$ around $p$. 
\begin{enumerate}
\item Suppose $p$ is reduced and non-degenerate singularity. 
Then we may assume that $\omega = \lambda y(1+o(1))\tn{d}x-x(1+o(1))\tn{d}y$ around $p$. 
Hence, 
\begin{align*}
& \textnormal{CS}(\sF,x=0,p) = \frac{1}{\lambda}, \textnormal{CS}(\sF,y=0,p) = \lambda \mbox{, and } \\
& \textnormal{Z}(\sF,x=0,p)=1=\textnormal{Z}(\sF,y=0,p).
\end{align*}
\item Suppose $p$ is a saddle-node. 
Then we may assume that 
\[\omega = y^{k+1}\tn{d}x-(x(1+\nu y^k)+yo(k))\tn{d}y\] 
around $p$. Hence,  $\textnormal{CS}(\sF,y=0,p) = 0$ and $\textnormal{Z}(\sF,y=0,p) = 1$.
Suppose there exists a weak separatrix, then $\textnormal{CS}(\sF,x=0,p) = \nu \mbox{ and } \textnormal{Z}(\sF,x=0,p) = k+1$.

\item Suppose $\sF$ has a positive rational eigenvalue $\lambda = \frac{m}{n}\in\bQ^+$ in the reduced form at $p$. 
Then we may assume $\omega = \lambda y(1+o(1))\tn{d}x-x(1+o(1))\tn{d}y$ around $p$. 
Hence, 
\begin{align*}
& \textnormal{CS}(\sF,x=0,p) = \frac{1}{\lambda},\,\textnormal{CS}(\sF,y=0,p) = \lambda,\,\mbox{and } \\
& \textnormal{Z}(\sF,x=0,p)=1=\textnormal{Z}(\sF,y=0,p).
\end{align*}
Moreover, for the invariant curve $C$ through $p=(0,0)$ which is locally defined by $ax^m(1+o(1))-y^n(1+o(1))=0$, we have 
\begin{align*}
& \textnormal{CS}(\sF,C,p) = mn \mbox{ and } \\
& \textnormal{Z}(\sF,C,p) = 1-(m-1)(n-1). 
\end{align*}

\end{enumerate}
\end{lem}
\begin{pf}
We will only show the case $(3)$. 
For the first two cases, the computations are similar. 
(See also \cite[page 30-31 in Chapter 3]{brunella2015birational}.) 

For the case $f=x$, note that $1\cdot\omega = (\lambda y(1+o(1)))\tn{d}x+x(-(1+o(1))\tn{d}y)$. 
Then we take $g=1$, $h=\lambda y(1+o(1))$, and $\eta = -(1+o(1))\tn{d}y$. 
Thus we have 
\[\textnormal{CS}(\sF,x=0,p) = \textnormal{res}_{p}\frac{-1}{h}\eta\vert_{f=0} 
= \textnormal{res}_{p}\frac{1}{\lambda y(1+o(1))}(1+o(1))\tn{d}y\vert_{x=0} = \frac{1}{\lambda}\]
and 
\begin{align*}
Z(\sF,x=0,p) &= \textnormal{vanishing order of } \frac{h}{g}\vert_{f=0} \\
&= \textnormal{vanishing order of } \frac{\lambda y(1+o(1))}{1}\vert_{x=0} \\
&= 1.
\end{align*}

For the case $f=y$, it is similar. 

For the case $C = (ax^m(1+o(1))-y^n(1+o(1))=0)$, we may assume, for simplicity, that $C = (f=ax^m-y^n=0)$ and $\omega = my\tn{d}x-nx\tn{d}y$. 
Note that 
\begin{align*}
ax^{m-1}\omega &= amx^{m-1}y\tn{d}x - anx^m\tn{d}y \\
&= amx^{m-1}y\tn{d}x - ny^n\tn{d}y - nf\tn{d}y \\
&= y\tn{d}f-nf\tn{d}y.
\end{align*}
So we take $g = ax^{m-1}$, $h=y$, and $\eta = -n\tn{d}y$. 
Thus, \[\textnormal{CS}(\sF,C,p) = \textnormal{res}_{p}\frac{-1}{h}\eta\vert_{f=0} 
= \textnormal{res}_{p}\,\frac{1}{y}n\tn{d}y\vert_{f=0} = mn\]
and 
\begin{align*}
\textnormal{Z}(\sF,C,p) &= \textnormal{vanishing order of } \frac{h}{g}\vert_{f=0} \\
&= \textnormal{vanishing order of } \frac{y}{ax^{m-1}}\vert_{f=0} \\
&= 1-(n-1)(m-1). 
\end{align*}
\qed
\end{pf}

\subsection{\texorpdfstring{$\ve$}{epsilon}-canonical singularities}
From Example~\ref{key_eg}, we know that there is no Seidenberg theorem for a family. 
Nevertheless, in \cite{pereira2019effective}, we can still have a Seidenberg-type theorem for a family by allowing $\ve$-canonical singularities, which will be defined as follows. 
\begin{defn}[{\cite[Definition 4.1]{pereira2019effective}}]
Let $(X,\sF)$ be a foliated surface with $X$ smooth and let $\pi: Y\rw X$ be a birational morphism. 
Denote $\sG$ be the pullback foliation of $\sF$ via $\pi$. 
If $E$ is a prime exceptional divisor of $\pi$, then we define the \emph{adjoint discrepancy} of $\sF$ along $E$ to be the function 
\[a(\sF,E) : [0,\infty)\rw \bR\] 
which sends $t$ to $\ord_E\big((1-t)(K_\sG-\pi^*K_\sF) + t(K_Y-\pi^*K_X)\big)$. 
\end{defn}

\begin{defn}[{\cite[Definition 4.2]{pereira2019effective}}]
A point $x$ on a foliated surface $(X,\sF)$ with $X$ smooth is called $\ve$-canonical if the adjoint discrepancy of $\sF$ along any divisor $E$ over $x$ satisfies $a(\sF,E)(t)\geq 0$ for every $t\geq \ve$. 
\end{defn}

\begin{lem}[{\cite[Lemma 7.3]{pereira2019effective}}]
Let $(\cX,\cF)$ be a smooth family over an algebraic variety $T$ such that all fibers are smooth projective surfaces. 
If $0<\ve<\frac{1}{4}$, then the subset of $T$ corresponding to foliations with isolated and $\ve$-canonical singularities is a Zariski open subset of $T$. 
\end{lem}

\begin{prop}[{\cite[Proposition 7.4]{pereira2019effective}}]\label{ve_Sei}
Given a smooth family $(\cX,\cF)$ over $T$ such that all fibers are smooth projective surfaces and a real number $\ve$ with $0<\ve < \frac{1}{4}$. 
Then there exists a non-empty Zariski open subset $U\subset T$ and a family $(\cY,\cG)$ over $U$ obtained from $\cF\vert_U$ by a finite composition of blowups over (multi)-sections such that, for all $t\in U$, the foliation $\cG_t$ has at worst $\ve$-canonical singularities. 
In particular, $\cY_t$ is smooth and $\cG_t$ has at worst log canonical singularities for all $t\in U$. 
\end{prop}

\subsection{Quot scheme}
In this section, we recall some properties of Quot schemes. 

\begin{defn}
Let $X \rw S$ be a scheme of finite type over a noetherian base scheme $S$ and $\cE$ be a coherent sheaf on $X$. 
We consider the functor 
\[\cQ uot_{\cE/X/S} : (\textnormal{Sch}/S)^{\textnormal{op}} \rw (\textnormal{Sets})\]
sending $T\rw S$ to the set 
\[\cQ uot_{\cE/X/S}(T) = \left\{(\cF,q)\Big\vert 
\begin{array}{l}
\cF\in\textnormal{Coh}(X_T), \supp{\cF} \textnormal{ is proper over } T \\
\cF \textnormal{ is flat over } T, q:\cE_T\rw \cF \textnormal{ is surjective.}
\end{array}\right\}\Big/\sim
\]
where $X_T = X\times_ST$, $\cE_T = \textnormal{pr}_X^*\cE$ under the projection $\textnormal{pr}_X : X_T \rw X$, and $(\cF , q ) \sim ( \cF', q')$ if $\ker q = \ker q'$. 
\end{defn}

\begin{defn}\label{Hilbert_poly}
Let $X \rw S$ be a scheme of a finite type over a noetherian base scheme $S$ and $\cA$ be a relatively ample line bundle. 
For any coherent sheaf $\cF$ on $X$ flat over $S$, we define the Hilbert polynomial $\Phi_\cF$ with respect to $\cA$ as \[\Phi_\cF(m) = \chi(\cF_s(m)) = \sum_{i=0}^{\dim X}(-1)^ih^i(X_s, \cF_s\otimes\cA_s^{\otimes m})\]
for any fixed $s\in S$. 
Note that this definition is independent of the choice of $s\in S$. 
\end{defn}

\begin{prop}
Fix any relatively ample line bundle $\cA$ on $X$. 
There is a natural stratification: 
\[\cQ uot_{\cE/X/S} = \coprod_{\Phi\in\bQ[\lambda]}\cQ uot_{\cE/X/S}^{\Phi,\cA}\]
where $\cQ uot_{\cE/X/S}^{\Phi,\cA}(T) = \{(\cF,q)\in\cQ uot_{\cE/X/S}(T)\,\vert\, \Phi_\cF = \Phi\}$. 
\end{prop}

\begin{thm}[Grothendieck]\label{quot_exist}
The functor $\cQ uot_{\cE/X/S}^{\Phi,\cA}$ is representable by a projective scheme $\textnormal{Quot}_{\cE/X/S}^{\Phi,\cA}$ over $S$.
\end{thm}

\section{Log canonical models} 
In light of the Example~\ref{key_eg}, we would like to allow some mild singularities in order to get local-closedness in a family. 
We introduce the following definition.

\begin{defn}\label{defn_gcm}
A foliated normal surface $(X,\sF)$ is called a \emph{log canonical model} if the following conditions hold:
\begin{enumerate}
\item $\sF$ has \emph{mild} log canonical foliation singularities, that is, $\sF$ is a foliation with only log canonical foliation singularities such that all points at which $\sF$ is \emph{not} canonical are either smooth points of $X$ or cyclic quotient singularities of $X$.
\item $K_\sF$ is nef. 
\item $K_\sF\cdot C = 0$ implies $C^2\geq 0$ for any irreducible curve $C$. 
\end{enumerate}
\end{defn}

In Lemma~\ref{pos_rat} and Remark~\ref{mlc}, we study more properties of mild log canonical foliation singularities. 
Actually, there are log canonical foliation singularities which are not mild. 
\begin{eg}[{\cite[Example I.2.5]{mcquillan2008canonical}}]
We consider the germ of the normal surface $(x^{n+2}+y^{n+1}+xz^n=0)$ at the origin for $n\geq 2$. 
It is not a rational singularity but the foliation given by 
\[\partial = n(n+1)x\p{x}+n(n+2)y\p{y}+(n+1)^2z\p{z}\]
is log canonical. 
\end{eg}

\begin{defn}
We say a foliated normal surface $(X,\sF)$ is a \emph{canonical model} if it is a log canonical model with only canonical foliation singularities. 
\end{defn}

\begin{lem}\label{ample}
Let $(X,\sF)$ be a log canonical model. 
If the foliation $\sF$ is of general type, then $K_{\sF}$ is numerically ample, that is $K_{\sF}^2>0$ and $K_{\sF}\cdot C>0$ for any irreducible curve $C$ on $X$. 
\end{lem} 
\begin{pf}
The proof follows closely from the proof of \cite[Lemma 1.11]{hacon2021birational}. \qed
\end{pf}

To describe the resolutions of canonical foliation singularities and mild log canonical foliation singularities, we recall the following definitions. 
(See also \cite[Definition 5.1 and 8.1]{brunella2015birational} and \cite[Definition III.0.2 and III.2.3]{mcquillan2008canonical}.)

\begin{defn}\label{-2curve}
Given a foliated surface $(X,\sF)$ with $X$ smooth. 
\begin{enumerate}
\item $E = \bigcup_{i=1}^sE_i$ is called a string if 
\begin{enumerate}
\item each $E_i$ is a smooth rational curve and 
\item $E_i\cdot E_j = 1$ if $|i-j|=1$ and $0$ if $|i-j|\geq 2$. 
\end{enumerate}
\item If, moreover, $E_i^2\leq -2$ for all $i$, then we call $E$ a Hirzebruch-Jung string. 
\item A curve $C\subset X$ is called $\sF$-exceptional if 
\begin{enumerate}
\item $C$ is a $(-1)$-curve and 
\item the contraction of $C$ introduces $(X',\sF')$ with only \emph{reduced} singularities. 
\end{enumerate}

\item An $\sF$-invariant curve $C$ is called a $(-1)$-$\sF$-curve (resp. $(-2)$-$\sF$-curve) if 
\begin{enumerate}
\item $C$ is a smooth rational curve and 
\item $\textnormal{Z}(\sF,C)=1$ (resp. $\textnormal{Z}(\sF,C)=2$). 
\end{enumerate}

\item We say $C=\bigcup_{i=1}^sC_i$ is an $\sF$-chain if 
\begin{enumerate}
\item $C$ is a Hirzebruch-Jung string, 
\item each $C_i$ is $\sF$-invariant, 
\item $\textnormal{Sing}{\sF}\cap C$ are all reduced and non-degenerate, and 
\item $\textnormal{Z}(\sF, C_1) = 1$ and $\textnormal{Z}(\sF,C_i) = 2$ for all $i\geq 2$. 
\end{enumerate}

\item If an $\sF$-invariant curve $E$ is not contained in an $\sF$-chain $C$ but meets the chain, 
then we call $E$ the tail of the chain $C$. 

\item $C$ is called a bad tail if 
\begin{enumerate}
\item $C$ is a smooth rational curve with $\textnormal{Z}(\sF,C) = 3$ and $C^2\leq -2$ and 
\item $C$ intersects two $(-1)$-$\sF$-curves whose self-intersections are both $-2$. 
\end{enumerate}
\end{enumerate}
\end{defn}

\begin{thm}[{\cite[Theorem III.3.2]{mcquillan2008canonical}}]\label{struct_can}
Given any canonical model $(X,\sF)$. 
Let $\pi: Y \rw X$ be the minimal resolution at points at which $\sF$ is singular. 
Let $\sG$ be the pullback foliation on $Y$. 
Then all connected components of exceptional divisors of $\pi$ are one of the following types:
\begin{enumerate}
\item ($A_n$ type) A string of smooth rational curves. 
More precisely, it consists of either a $\sG$-chain, two $(-1)$-$\sG$-curves of self-intersections $-2$ joined by a bad tail, or a chain of $(-2)$-$\sG$-curves.
\item ($D_n$ type) Two $(-1)$-$\sG$-curves of self-intersections $-2$ joined by a bad tail which itself connects to a chain of $(-2)$-$\sG$-curves. 
\item (e.g.l.) Elliptic Gorenstein leaves. That is, either a cycle of $(-2)$-$\sG$-curves or a rational curve with only one node.
\end{enumerate}
\end{thm}

Now we first study mild log canonical foliation singularities on \emph{smooth} points.
\begin{lem}\label{pos_rat}
Let $(X,\sF,p)$ be a germ of a foliated surface with $p\in X$ a smooth point. 
Suppose $\sF$ at $p$ has the positive rational eigenvalue. 
Then there is a sequence of blowups $\pi : (Y,\sG) \rw (X,\sF,p)$ with $\sG$ having only reduced singularities on the exceptional divisor $E$. 
Moreover we have the following:
\begin{enumerate}
\item The support of the exceptional divisor $\bigcup_{i=1}^sE_i$ is a string, 
\item All $E_i$'s except for one $j\in\{1,\ldots,s\}$ are invariant; while this $E_j$ is non-invariant with $\textnormal{tang}(\sG,E_j)=0$. 
\item The string $\bigcup_{i=1}^{j-1}E_i$ (resp. $\bigcup_{i=j+1}^{s}E_i$) is a $\sG$-chain with the first curve $E_{j-1}$ (resp. $E_{j+1}$). 
\item $\pi^*K_\sF = K_\sG + E_j$. 
\end{enumerate}
\end{lem}
\begin{pf}
Since the eigenvalue is defined up to reciprocal, we may assume that the eigenvalue is $\frac{m}{n}\in\bQ^+$ with $\textnormal{gcd}(m,n)=1$ and $m\geq n$. 
We proceed by induction on the pair $(m,n)\in\bN\times\bN$ (by lexicographic ordering). 
\begin{enumerate}
\item We may assume that $p=(0,0)$ and that $\sF$ is generated by 
\[v = nx\p{x}+my\p{y}\] 
around $p$. 
Let $\varphi : (Y,\sG)\rw (X,\sF)$ be the blowup at $p$ and $E_0$ be the exceptional divisor. 
Notice that there are two singularities on the exceptional curve locally generated by 
\[v_1 = nx'\p{x'}+(m-n)y'\p{y'}\]
and 
\[v_2 = (n-m)x'\p{x'}+my'\p{y'}.\]

\item For the base case $(m,n) = (1,1)$, by the computation above, the exceptional curve $E_0$ has no foliation singularity and is not invariant with $\textnormal{tang}(\sG,E_0) = 0$. 
Moreover, we have $\pi^*K_\sF = K_\sG+E_0$. 
Thus, we complete the proof of the base case. 

\item Assume $(m,n)\neq (1,1)$ and all pairs (strictly) less than $(m,n)$ have the desired properties. 
Then $E_0$ is invariant and $\pi^*K_\sF = K_\sG$. 
Notice that $v_1$ still has positive rational eigenvalue with corresponding pair $(m-n,n)$ or $(n,m-n)$ while $v_2$ is reduced. 

By induction hypothesis, $v_1$ has a sequence of blowups $\psi : (Z,\sH) \rw (Y,\sG)$ with exceptional divisor $E = \bigcup_{i=1}^sE_i$. 
Since the proper transform $\wt{E}_0$ of $E_0$ on $Z$ meets $E$, $\wt{E}_0$ can only meet $E$ as a tail of an $\sH$-chain. 
Note that $\wt{E}_0^2\leq E_0^2-1=-2$. 
Thus, $\wt{E}_0$ is a member of that $\sH$-chain. 
Also $\psi^*\pi^*K_\sF = \psi^*K_\sG = K_\sH+E_i$ for some $i$.
\qed
\end{enumerate}
\end{pf}

\begin{rmk}\label{mlc}
By \cite[Fact I.2.4 and its remarks]{mcquillan2008canonical}, a mild log canonical foliation singularity $p$ (see Definition~\ref{defn_gcm}) is formally locally of the form $\partial = x\p{x}+\lambda y\p{y}$ quotient by $\bZ/n$ for some $n\in\bZ$ and $\lambda\in\bQ^+$. 
Also $K_\sF$ is Cartier around $p$. 
Let $\pi : Y\rw X$ be the \emph{minimal} resolution of $X$ at $p$ and $\sG$ be the pullback foliation on $Y$. 
Then, after some computations, we have the following:
\begin{enumerate}
\item The support of the exceptional divisor $E = \bigcup_{i=1}^sE_i$ is a string. 
\item For general $\lambda\in\bQ^+$, all $E_i$'s are $(-2)$-$\sG$-curves, exactly one foliation singularity $q$ on $E$ is not reduced, and $\sG$ has the positive rational eigenvalue at $q$.  
\item Fix the action of $\bZ/n$ on $\partial$, for any $j\in\{1,\ldots,s\}$, there is a $\lambda\in\bQ^+$ such that $E_j$ is $\sG$-non-invariant with $\textnormal{tang}(\sG,E_j)=0$, the string $\bigcup_{i=1}^{j-1}E_i$ (resp. $\bigcup_{i=j+1}^{s}E_i$) is a $\sG$-chain with the first curve $E_{j-1}$ (resp. $E_{j+1}$), and $\pi^*K_\sF = K_\sG + E_j$. 
\end{enumerate}
\end{rmk}

\begin{eg}
Consider the foliation of the form $\partial = x\p{x}+\lambda y\p{y}$ on the germ $(\bC^2/\mu_7,0)$ of the type $\frac{1}{7}(1,4)$. 
Let $Y\rw X$ be the minimal resolution of singularities. 
Then we have two irreducible exceptional divisors $E_1$ and $E_2$, and we need three copies of $\bC^2$ to cover $Y$. 
More explicitly, $Y = \bC^2_{\xi_0,\eta_0}\cup\bC^2_{\xi_1,\eta_1}\cup\bC^2_{\xi_2,\eta_2} = \bA^2_0\cup\bA^2_1\cup\bA^2_2$ 
where \[\begin{array}{lll}
\xi_0=x^{-4}y, & \xi_1 = x^{-1}y^2, & \xi_2 = y^7 \\
\eta_0 = x^7, & \eta_1 = x^4y^{-1}, & \eta_2 = xy^{-2}
\end{array}\]
Thus, on $\bC^2_{\xi_0,\eta_0}$, we have 
\[\partial_0 = (\lambda-4)\xi_0\p{\xi_0}+7\eta_0\p{\eta_0}.\]
On $\bC^2_{\xi_1,\eta_1}$, we have 
\[\partial_1 = (2\lambda-1)\xi_1\p{\xi_1}+(4-\lambda)\eta_1\p{\eta_1}.\]
And on $\bC^2_{\xi_2,\eta_2}$, we have 
\[\partial_2 = 7\lambda\xi_2\p{\xi_2}+(1-2\lambda)\eta_2\p{\eta_2}.\]
Assume that $E_1$ is on $\bA^2_0\cup\bA^2_1$ and thus defined by $\eta_0=0$ on $\bA^2_0$ and by $\xi_1=0$ on $\bA^2_1$. 
So $E_2$ is on $\bA^2_1\cup\bA^2_2$ and defined by $\eta_1=0$ on $\bA^2_1$ and by $\xi_2=0$ on $\bA^2_2$. 

Therefore, if $\lambda = 4$, then $E_1$ is non-invariant. 
Also if $\lambda = \frac{1}{2}$, then $E_2$ is non-invariant. 
For other $\lambda\in\bQ^+$, the exceptional curves $E_1$ and $E_2$ are both invariant, and precisely one of three foliation singularities has positive rational eigenvalue for the foliation. 
\end{eg}

Now we summarize the descriptions of the singularities on a log canonical model. 
\begin{thm}\label{description_GenCanModel}
Let $(X,\sF)$ be a log canonical model. 
Let $\pi: Y \rw X$ be the minimal resolution at points at which $\sF$ is singular. 
Let $\sG$ be the pullback foliation on $Y$. 
Write $K_\sG+D = \pi^*K_\sF+E$ where $D$ and $E$ are effective and have no common components. 
Then $\sG$ has at mild log canonical singularities, $D$ is reduced, the support of $E$ is the union of $\sG$-chains which are disjoint from $D$, and all irreducible components of $D$ are pairwise disjoint with tangency order zero. 
Moreover, the connected components of $\pi$-exceptional divisors belong to one of the following types: 
\begin{enumerate}
\item ($A_n$ type) A chain of rational invariant curves consisting of an $\sG$-chain, two $(-1)$-$\sG$-curves of self-intersections $-2$ joined by a bad tail, a chain of $(-2)$-$\sG$-curves on which there is at most one non-reduced singularities, or a non-invariant smooth rational curve with tangency order zero.
\item ($D_n$ type) Two $(-1)$-$\sG$-curves of self-intersections $-2$ joined by a bad tail which itself connects to a chain of $(-2)$-$\sG$-curves on which there are only reduced singularities. 
\item (e.g.l.) Elliptic Gorenstein leaves. That is, either a cycle of $(-2)$-$\sG$-curves or a rational curve with only one node. 
\end{enumerate}
\end{thm}

\begin{rmk}
Only the eigenvalues of foliation singularities on the chain of $(-2)$-$\sG$-curves can vary in a family. 
\end{rmk}

\section{Boundedness of log canonical models}
\subsection{Local contributions of log canonical models}
From now on, we assume that the foliation $\sF$ is of general type, that is, $K_\sF$ is big. 
Fix $\chi(mK_\sF) = P(m)$ where $P : \bZ_{\geq 0} \rw \bZ$ is an integral valued function. 
As in \cite{hacon2021birational}, we would like to study the following singular Riemann-Roch theorem and get some boundedness results. 

\begin{thm}[{\cite{reid1987canonical},\cite[Section 3]{langer2000chern}}]\label{RR}
Let $X$ be a normal complete surface and $D$ be a Weil divisor on $X$. 
Then we have 
\[\chi(D) = \frac{1}{2}D\cdot(D-K_X) + \chi(\cO_X)+\sum_{x\in\textnormal{Sing}(X)} a(x,D)\]
where $a(x,D)$ is a local contribution of $\cO_X(D)$ at $x$ depending only on the local isomorphism class of the reflexive sheaf $\cO_X(D)$ at $x$.
\end{thm}

If $X$ has only quotient singularities, then the theorem follows from \cite{reid1987canonical}. 
In general, the theorem follows from \cite[Section 3]{langer2000chern}. 

\begin{prop}\label{local_contribution}
Let $m$ be a non-negative integer. 
\begin{enumerate}
\item If $D$ is Cartier at $x$, then $a(x,D)=0$. 
\item If $(X,\sF,x)$ is a mild log canonical non-canonical foliation singularity, then 
\[a(x,mK_\sF) = 0 \mbox{ for all } m.\]
\item If $(X,\sF,x)$ is a canonical foliation singularity, which is $\bQ$-Gorentein but it is not terminal. 
Then we have two cases:
\begin{enumerate}
\item $\sF$ is Gorenstein and $a(x,mK_\sF) = 0$ for all $m$.
\item $\sF$ is $2$-Gorenstein and 
\[a(x,mK_\sF) = \left\{\begin{array}{ll}
0 & \mbox{ for $m$ even} \\
\frac{-1}{2} & \mbox{ for $m$ odd.}
\end{array}\right.\] 
\end{enumerate}
\item If $(X,\sF,x)$ is a canonical foliation singularity and $\sF$ is not $\bQ$-Gorenstein at $x$, then we have 
\[a(x,mK_\sF) = \left\{\begin{array}{ll}
0 & \mbox{ for $m=0$} \\
-1 & \mbox{ for positive $m$.}
\end{array}\right.\] 

\item If $(X,\sF,x)$ is a terminal foliation singularity, then $a(x,K_\sF) = \frac{1-n}{2n}$ where the surface singularity at $x$ is of the type $\frac{1}{n}(1,q)$ for some $n$ and $q$. 
\end{enumerate}
\end{prop}
\begin{pf}
(1) follows from the definition. 
(2) follows from (1) since  $K_\sF$ is Cartier at these foliation singularities. 
(3) is \cite[Proposition 2.4]{hacon2021birational} and (4) is \cite[Proposition 2.6]{hacon2021birational}. 
(5) follows from \cite[Corollary 2.2]{hacon2021birational} and \cite[Example 5.6]{langer2000chern}. 
\qed
\end{pf}

\begin{rmk}
In addition to any singularities on canonical models, we only allow \emph{mild} log canonical foliation singularities, which have no local contributions by Proposition~\ref{local_contribution}. 
Thus, we have the following similar proposition as \cite[Proposition 4.1]{hacon2021birational} for log canonical models. 
\end{rmk}

\begin{prop}\label{bdd}
Fix a function $P : \bZ_{\geq 0} \rw \bZ$. 
There exist constants $B_1$, $B_2$, $B_3$, and $B_4$ such that if $(X,\sF)$ is a log canonical model with Hilbert function $\chi(mK_\sF) = P(m)$ for all $m\in\bZ_{\geq 0}$, then $K_\sF^2 = B_1$, $K_\sF\cdot K_X = B_2$, $\chi(\cO_X) = B_3$, and the number of cusp singularities is equal to $B_4$. 
Moreover, there exists some constants $C_1$ and $C_2$ such that the number of terminal and dihedral singularities is at most $C_1$ and the index at any terminal foliation singularity is at most $C_2$. 
\end{prop}
\begin{pf}
The proof follows closely from the proof of \cite[Proposition 4.1]{hacon2021birational}. 
For the reader's convenience, we include the proof. 

By Theorem~\ref{RR}, we have 
\[P(m) = \chi(mK_\sF) = \frac{1}{2}mK_\sF\cdot (mK_\sF-K_X)+\chi(\cO_X)+\sum_{x\in\textnormal{Sing}(X)} a(x,mK_\sF)\]
Since $K_\sF^2$ is the quadratic term of $P$ and $K_\sF\cdot K_X$ is the linear term of $P$, both of them are fixed. 
Also $\chi(\cO_X) = P(0)$ is also fixed. 
If we put $m=1$, then we get the number 
\[\sum_{x\in\textnormal{Sing}(X)} a(x,K_\sF)\]
is fixed. 

Let $\Sigma_1$ be the set of terminal foliation singularities, $\Sigma_2$ be the set of dihedral quotient singularities, and $\Sigma_3$ be the set of cusp singularities. 
Also let $\Sigma = \Sigma_1\cup\Sigma_2\cup\Sigma_3$. 
Then by Proposition~\ref{local_contribution}, we have  
\[-\sum_{x\in\textnormal{Sing}(X)} a(x,K_\sF) = \sum_{x\in\Sigma_1}\frac{n_x-1}{2n_x} + \sum_{x\in\Sigma_2}\frac{1}{2}+\sum_{x\in\Sigma_3}1 \geq \frac{1}{4}\vert\Sigma\vert.\] 
So the number of singularities is bounded and hence 
\[\sum_{x\in\Sigma_1}\frac{1}{n_x} = \vert\Sigma\vert+\vert\Sigma_3\vert+2\sum_{x\in\textnormal{Sing}(X)} a(x,K_\sF)\]
assumes finitely many values. 
By \cite[Lemma 3.4]{hacon2021birational}, the indices at the terminal foliation singularities are bounded by $C_2$. 

Finally, let $m = (2C_2)!$. 
Then by Theorem~\ref{RR} and Proposition~\ref{local_contribution}, we have 
\[P(m) = \chi(mK_\sF) = \frac{1}{2}mK_\sF\cdot(mK_\sF-K_X)+\chi(\cO_X) - \vert\Sigma_3\vert.\]
Hence, the number of cusps is fixed. 
\qed
\end{pf}

\subsection{Minimal partial du Val resolution}
Fix a function $P : \bZ_{\geq 0} \rw \bZ$. 
Given any log canonical model $(X,\sF)$ of general type with $\chi(mK_\sF) = P(m)$ for all $m\in\bZ_{\geq 0}$. 
We first consider the minimal resolution $g : X^{\tn{m}} \rw X$ at points at which $\sF$ is \emph{not} terminal. 
Let $\sF^{\tn{m}}$ be the pullback foliation on $X^{\tn{m}}$. 
Then let $h : X^{\tn{m}} \rw Y$ be the (traditional) relative canonical model of $K_{X^{\tn{m}}}$ over $X$ which is given by a sequence of contractions of smooth rational curves $C$ in the fiber of $g$ with $C^2 = -2$. 
Put $\sG$ be the pushforward foliation on $Y$. 
Thus we have $f : (Y,\sG) \rw (X,\sF)$ a morphism of foliated surfaces such that $g = f\circ h$. 

By Theorem~\ref{description_GenCanModel}, we have $K_{\sF^{\tn{m}}}+D^{\tn{m}} = g^*K_\sF+E$ where each irreducible component of $E$ is an $(-1)$-$\sF^{\tn{m}}$-curve with self-intersection $-2$ and $D^{\tn{m}} = \sum_i D^{\tn{m}}_i$ is an effective reduced divisor such that all irreducible components $D^{\tn{m}}_i$'s are pairwise disjoint and non-invariant with $\textnormal{tang}(\sF^{\tn{m}},D^{\tn{m}}_i)=0$. (See also Remark~\ref{mlc}.)
Moreover, $E$ is contracted by $h$ since each irreducible component is $(-1)$-$\sF^{\tn{m}}$-curve with self-intersection $-2$ and any two $(-1)$-$\sF^{\tn{m}}$-curves are disjoint. 
Thus, we have $K_\sG+D = f^*K_\sF$ where $D = h_*D^{\tn{m}}$. 
\[\xymatrix{(X^{\tn{m}},\sF^{\tn{m}}) \ar[rd]^-h \ar[dd]_-g & \\
 & (Y,\sG) \ar[ld]^-f \\
(X,\sF) &}\]

\begin{defn}\label{mpdvr}
We call $f : (Y,\sG,D) \rw (X,\sF)$ as above the \emph{minimal partial du Val resolution} of the given log canonical model where $(Y,\sG,D)$ is a foliated triple (see Definition~\ref{fol_triple}). 
\end{defn}

\begin{rmk}\label{unique}
If $(Y,\sG,D)$ is the minimal partial du Val resolution of some log canonical models $(X,\sF)$ of general type, then $(X,\sF)$ is unique. 
More precisely, if $(Y,\sG,D)$ is the minimal partial du Val resolution of two log canonical models $(X_1,\sF_1)$ and $(X_2,\sF_2)$. 
We have $f_1^*K_{\sF_1} = K_\sG+D = f_2^*K_{\sF_2}$. 
Let $C$ be a curve contracted by $f_1$. 
So we have $f_2^*K_{\sF_2}\cdot C = f_1^*K_{\sF_1}\cdot C=0$. 
Since $K_{\sF_2}$ is numerically ample by Lemma~\ref{ample}, $C$ is also contracted by $f_2$. 
\end{rmk}

\begin{defn}\label{set_MPdVR}
Fix a function $P : \bZ_{\geq 0} \rw \bZ$. 
Let $\cS_P$ be the set of foliated triples $(Y,\sG,D)$ which are minimal partial du Val resolutions of log canonical models $(X,\sF)$ of general type with $\chi(mK_\sF) = P(m)$. 
Also let $\cS_P^C\subset \cS_P$ be the subset consisting of foliated triples $(Y,\sG,\emptyset)$ which are minimal partial du Val resolutions of canonical models $(X,\sF)$ of general type with $\chi(mK_\sF)=P(m)$. 
\end{defn}

\begin{defn}\label{family_fol}
Let $\mu: \cY \rw T$ be a projective morphism between two quasi-projective varieties of finite type. We say a foliation $\cF$ (of rank one) on $\cY$ is a family of foliations over $T$ if $\Omega_{\cY}^{[**]} \rw \Omega_\cF$ factors through the relative hull $\Omega_{\cY/T}^{[**]}$ over $T$. 
Here the (relative) hull is defined in \cite[Definition 14 and 17]{kollar2008hulls}.
\end{defn}

Now we would like to prove the following boundedness for the set $\cS_P$. 
\begin{thm}\label{key_thm}
The set $\cS_P$ is bounded. 
More precisely, there exists a projective morphism $\mu : \cY \rw T$ and a family of foliations $\cF$ on $\cY$ over $T$ where $\cY$ and $T$ are quasi-projective varieties of finite type such that, for any $(Y,\sG,D)\in\cS_P$, there is a $t\in T$ and an isomorphism $\rho : Y\rw \cY_t$ such that $\sG(-D) \cong \rho^*(\cF\vert_{\cY_t})$.
\end{thm}

\begin{pf}
The proof follows closely from the proof of \cite[Theorem 0.1]{chen2021boundedness} for canonical models. 
For the reader's convenience, we include the proof. 
The difference is that now we possibly have non-invariant exceptional divisors $D = \sum D_j$. 

Let $(X,\sF)$ be the log canonical model of $(Y,\sG,D)$ (see Remark~\ref{unique}) with the associated morphism $f : (Y,\sG,D) \rw (X,\sF)$. 
Notice that $f^*K_\sF = K_\sG+D$. 
Let $E = \cup E_i$ be the exceptional divisor of $f : Y \rw X$. 
We divide the proof into several steps.

\begin{enumerate}
\item\label{adjoint_ample} \textit{Claim.} $K_Y+mi(\sG)(K_\sG+D)$ is nef (resp. ample) when $m \geq 3$ (resp. $m\geq 4$).
\begin{pf}
Given any irreducible curve $C$ on $Y$. 
\begin{enumerate}
\item If $C=E_i$ for some $i$, then 
\[(K_Y+mi(\sG)(K_\sG+D))\cdot C = K_Y\cdot C > 0\]
since $K_Y$ is ample over $X$. 
\item If $C$ is not contained in the exceptional divisor, then we have two cases: 
\begin{enumerate}
\item If $K_Y\cdot C\geq 0$, then 
\begin{align*}
(K_Y+mi(\sG)(K_\sG+D))\cdot C &\geq mi(\sG)(K_\sG+D)\cdot C \\
&= mi(\sG)K_{\sF}\cdot f_*C > 0 
\end{align*}
since $K_\sF$ is numerically ample (see Lemma~\ref{ample}). 

\item If $K_Y\cdot C<0$, by \cite[Proposition 3.8]{fujino2012minimal}, every $K_Y$-negative extremal ray is spanned by a rational curve $C$ with $0<-K_Y\cdot C\leq 3$. 
Thus, if we write $C \equiv \sum a_iC_i+B$ where $B\cdot K_Y\geq 0$, the $C_i$'s are $K_Y$-negative extremal rays, the $a_i$'s are non-negative, and at least one of $a_i$ is positive. 
By what we have seen above, we have 
\[(K_Y+mi(\sG)(K_\sG+D))\cdot B\geq 0.\]
Thus, we have 
\begin{align*}
(K_Y+mi(\sG)(K_\sG+D))\cdot C &\geq \sum_ia_i \Big(K_Y\cdot C_i+mi(\sG)(K_\sG+D)\cdot C_i\Big) \\
&\geq \sum_ia_i(- 3+m) \geq 0 \mbox{ (resp. $>0$)}.
\end{align*}
\end{enumerate} 
\end{enumerate}
When $m\geq 4$, we notice that 
\[K_Y+mi(\sG)(K_\sG+D) = K_Y+3i(\sG)(K_\sG+D) + (m-3)i(\sG)(K_\sG+D)\] 
is big since it is the sum of a nef divisor and a big divisor. 
\end{pf}

\item \textit{Claim.} $i(K_Y)$ and $i(\sG)$ are bounded at foliation singularities. 
\begin{pf}
By Proposition~\ref{bdd}, $i(K_Y)$ and $i(\sG)$ are bounded in terms of $\chi(mK_\sF)$ at terminal foliation singularities. Other foliation singularities are at du Val singularities at which $\sG$ is canonical. Then we have $i(K_Y) = 1$ and $i(\sG) \leq 2$ by Proposition~\ref{local_contribution}(3).
\end{pf}

\item \textit{Claim.} The embedding dimension at each point at which $\sG$ is singular is bounded. 
\begin{pf}
By Proposition~\ref{bdd}, the index at each point at which $\sG$ is terminal is bounded in terms of $\chi(mK_\sF)$. 
Actually, it is shown that the order of each cyclic group acting at each point at which $\sG$ is terminal is bounded. 
Thus, the embedding dimension at each point at which $\sG$ is terminal is bounded. 
See for example \cite[Corollary 2.5]{reid2012surface}. 
If the action is $\frac{1}{n}(1,q)$ and $\frac{n}{n-q} = [a_1,\ldots, a_k]$, then the embedding dimension is at most $k+2$. 

Other singularities on $Y$ are du Val singularities, at which the embedding dimensions are $3$. 
\end{pf}
Thus, we have a uniform bound $\delta$ such that $\delta_\zeta\leq\delta$ for $|\zeta|=2$. 
See Remark~\ref{delta}.

\item We consider the $\bQ$-Cartier divisor 
\begin{align*}
L &= \alpha i(K_Y)(K_Y+4i(\sG)(K_\sG+D)) - K_Y\\
&= (\alpha i(K_Y)-1)(K_Y+3i(\sG)(K_\sG+D))+(\alpha i(K_Y)+3)i(\sG)(K_\sG+D).
\end{align*}

Since $i(\sG)(K_\sG+D)$ is Cartier, 
we can fix an $\alpha$, which only depends on the Hilbert function $\chi(mK_\sF)$, such that for any irreducible curve $C$, we have 
\begin{align*}
L^2 &\geq (\alpha i(K_Y)+3)^2> \delta \mbox{ and}\\
L\cdot C &\geq \alpha i(K_Y)+3>\delta/2. 
\end{align*}
By Theorem~\ref{va_eff}, $K_Y+\lceil L\rceil = \alpha i(K_Y)(K_Y+ 4i(\sG)(K_\sG+D))$ separates points and tangents. 
Hence, by Theorem~\ref{sep}, it is very ample.

\item \textit{Claim.} The set of surfaces $X$ where $(X, \sF) \in \cS_P$ forms a bounded family.
\begin{pf}
Let $\varphi : Y \rw \bP^N$ be an embedding defined via $|K_Y + L|$. We may assume that $N$ is minimal in the sense that $\varphi(Y)$ is not contained in any hyperplane of $\bP^N$. 
In this case, it is known (for a reference, see \cite[Proposition 0]{eisenbud1987varieties}) that 
\[\deg\varphi(Y) \geq 1 + (N-\dim \varphi(Y)) = N-1.\]

By \cite[Lemma 3.5]{hacon2021birational}, we know that $\gamma K_{\sF}-K_{X}$ is pseudo-effective for 
\[\gamma = \max\left\{\frac{2K_{\sF}\cdot K_{X}}{K_{\sF}^2}+3i_\bQ(\sF),0\right\}\] 
which is bounded in terms of $\chi(mK_\sF)$. 
Thus, 
\begin{align*}
\deg\varphi(Y) &= \vol(K_Y+L) \\
&= \vol\Big(\alpha i(K_Y)(K_Y+ 4i(\sG)(K_\sG+D))\Big) \\
&\leq \vol\Big(\alpha i(K_Y)(K_{X}+ 4i(\sG)K_{\sF})\Big) \\
&\leq \vol\Big(\alpha i(K_Y)(K_{X}+ 4i(\sG)K_{\sF}) + \alpha i(K_Y)(\gamma K_{\sF}-K_{X}) \Big) \\
&= \vol\Big(\alpha i(K_Y)(\gamma + 4i(\sG))K_{\sF} \Big) \\
&= \Big(\alpha i(K_Y)(\gamma + 4i(\sG))K_{\sF} \Big)^2
\end{align*}
which is bounded in terms of $\chi(mK_\sF)$. 
Therefore $N$ is also bounded. 

Hence, by the boundedness result about Chow varieties, we have the boundedness of surfaces $Y$ where $(Y,\sG,D)\in\cS_P$. 
\end{pf}

\item Now, we would like to show the boundedness of $\sG(-D)$. 

We have shown that there exists a projective morphism $\mu : \cY\rw T$ where $\cY$ and $T$ are quasi-projective varieties of finite type such that for any $(Y,\sG,D)\in\cS_P$, there is a $t\in T$ such that $\cY_t\cong Y$ and by construction we have 
\[\xymatrix{\cY \ar@{^{(}->}[r] \ar[d]_-\mu & \bP^N_T \ar[ld] \\ T & }\]
and $\cO_Y(m(K_\sG+D) + nK_Y) = \cO_{\cY_t}(1)$ where $m=4\alpha i(K_Y)i(\sG)$ and $n = \alpha i(K_Y)$. 
Note that $\alpha$, $i(K_Y)$, and $i(\sG)$ are fixed and determined by the function $P$. 
Without loss of generality, we may assume that $T$ is the closure of $T'$ where 
\[T' = \{t\,\vert\, \cY_t\cong Y \mbox{ for some } (Y, \sG, D)\in \cS_P\}.\] 

\item We may assume that $T$ is reduced. 
\begin{claim}
We may assume that $T$ is irreducible. 
\end{claim}
\begin{pf}
Let $T = \bigcup_iT_i$ be the decomposition of $T$ into irreducible components. 
If there is a family of pre-foliations $\cF_i$ on each $\cY_i = \cY\times_T T_i$ such that for any $(Y,\sG,D)\in\cS_P$, there exist an $i$ and $t\in T_i$ with $(Y,\sG(-D))\cong (\cY_i,\cF_i)_t$, then the disjoint union of these families gives the required family. 
\end{pf}

\begin{claim}
We may assume that $T$ is smooth. 
\end{claim}
\begin{pf}
This is achieved by the existence of a finite stratification of $T$ by smooth locally closed subsets. 
\end{pf}

\item Recall that, by Serre's criterion, $\cY$ is normal if and only if it is $R_1$ and $S_2$. 
\begin{claim}
We may assume that $\cY$ is normal. 
\end{claim}
\begin{pf}
Suppose that $\cY_t$ is $R_1$ (resp. $S_2$), then there exists an open subset $U\subset T$ containing $t$ 
such that $\cY_U := \cY\times_T U$ is $R_1$ (resp. $S_2$) and every fiber $\cY_t$ is also $R_1$ (resp. $S_2$) for $t\in U$. 
(For a reference, see \cite[Th\'{e}or\`{e}me 12.2.4]{grothendieck1966elements}.) 
It follows that there is an open subset $V \subset T$ such that $\cY_V := \cY\times_T V$ is normal, every fiber $\cY_t$ is normal for $t\in V$, and for all $Y\in \cS$, there exists a $t\in V$ such that $Y\cong \cY_t$. 
Hence, by shrinking $T$, we may assume that $\cY$ is normal. 
\end{pf}

\item \textit{Claim.} We may assume that $K_\cY$ is $\bQ$-Cartierness. 
\begin{pf}
This is well-known; see, for example, \cite[Lemma 1.27]{simpson1994moduli}. 
We include proof for the convenience of the reader. 

Since $\cY$ is normal, the canonical sheaf $K_\cY$ is defined as a divisorial sheaf. 
By assumption, $i(K_Y)K_{\cY_t}$ is Cartier for 
\[t\in T^\circ := \{t\in T\,\vert\, \mbox{there exists } (Y,\sG,D)\in\cS_P \mbox{ such that } \cY_t\cong Y\}.\]
Note that $T^\circ$ is a dense subset of $T$. 
For any $t\in T^\circ$, we fix a $d$ such that $i(K_Y)K_{\cY_t}(d)$ is generated by global sections. 
Then by replacing $T$ by an open subset, we may assume that the map 
\[\mu_*\cO_\cY(i(K_Y)K_\cY(d))\otimes k(t) \rw H^0(\cY_t,i(K_Y)K_{\cY_t}(d))\]
is surjective by \cite[Theorem 12.8 and Corollary 12.9]{hartshorne1977algebraic}. 
Since $i(K_Y)K_{\cY_t}(d)$ is locally generated by some section of $H^0(\cY_t,i(K_Y)K_{\cY_t}(d))$, by lifting this section to a section of $H^0(\cY,i(K_\cY) K_\cY(d))$ using the surjection above and by Nakayama's lemma, $i(K_Y)K_\cY(d)$ is locally generated by one section on a neighborhood of $\cY_t$. 
Therefore, $i(K_Y)K_\cY$ is Cartier on a neighborhood of $\cY_t$. 
Shrinking $T$, we may assume that $i(K_Y)K_\cY$ is Cartier. 
\end{pf}

\item Now we consider the function 
\begin{align*}
\Phi(\ell) = {} & \chi(Y,\cO_X(K_\sG+D+\ell(m(K_\sG+D)+nK_Y))) \\
= {} & \chi(X,\cO_X((m\ell+1)(K_\sG+D)+n\ell K_Y)) \\
= {} & \frac{1}{2}\Big((m\ell+1)(K_\sG+D)+n\ell K_Y\Big)\cdot\Big((m\ell+1)(K_\sG+D)+(n\ell-1)K_Y\Big) \\
& + \chi(\cO_Y) + \sum_{y\in\textnormal{Sing}(Y)}a\big(y,(m\ell+1)(K_\sF+D)+n\ell K_Y\big) 
\end{align*}
Note that $m(K_\sG+D)$ and $nK_Y$ are Cartier, we have that 
\[a\big(y,(m\ell+1)(K_\sG+D)+n\ell K_Y\big) = a(y,K_\sG+D)\] 
for any singularity $y\in Y$. 
Also we have
\begin{align*}
\chi(\cO_Y) + \sum_{y\in\textnormal{Sing}(Y)}a(y,K_\sG+D) 
&= \chi(K_\sG+D) - \frac{1}{2}(K_\sG+D)\cdot(K_\sG+D-K_Y) \\
&= P(1) - \frac{1}{2}(K_\sG+D)\cdot(K_\sG+D-K_Y).
\end{align*}
Thus, we have  
\begin{align*}
\Phi(\ell) = {} & \frac{1}{2}\Big((m\ell+1)(K_\sG+D)+n\ell K_Y\Big)\cdot\Big((m\ell+1)(K_\sG+D)+(n\ell-1)K_Y\Big) \\
& + P(1) - \frac{1}{2}(K_\sG+D)\cdot(K_\sG+D-K_Y)
\end{align*}
is a polynomial with coefficients in terms of $P(1)$, $(K_\sG+D)\cdot K_Y = K_\sF\cdot K_X$, $(K_\sG+D)^2=K_\sF^2$, and $K_Y^2$. 
Since the function $P$ is fixed, we have $K_\sF^2$ and $K_\sF\cdot K_X$ are fixed by Proposition~\ref{bdd}. 
Moreover, $K_Y^2$ has only finitely many possibilities since the surfaces $Y$ are in a bounded family. 
Therefore, we have finitely many possible polynomials $\Phi(\ell)$. 

Notice that $\frac{m+1}{n-1}>\frac{m}{n} = 4i(\sG)\geq 4$. 
So we have $(m+1)(K_\sG+D)+(n-1)K_Y$ is ample by step~(\ref{adjoint_ample}) and thus $\Phi(1)=h^0(Y,\cO_Y((m+1)(K_\sG+D)+nK_Y))$. 

Let $\varphi$ be the maximum of $\Phi(1)$ among all possible polynomials $\Phi(\ell)$. 

\item Applying \cite[Theorem 21]{kollar2008hulls} to the projective morphism $\mu : \cY \rw T$ and the coherent sheaf $\Omega_{\cY/T}$, after shrinking $T$ further, we may assume that $\Omega_{\cY/T}$ has a relative hull $\Omega_{\cY/T}^{[**]} = \Omega_{\cY/T}^{**}$. 
Note that the relative hull behaves well under base change. 

For any polynomial $\Psi(\ell) := \Phi(\ell)-s$ where $0\leq s\leq \varphi$, we consider the Quot scheme 
\[Q := \textnormal{Quot}_{\Omega_{\cY/T}^{**}/\cY/T}^{\Psi(\ell),\cO_\cY(1)},\] 
which is projective over $T$. 
Let $\cU := \cY\times_T Q$. 
Then by a flattening stratification on $Q$, we may assume that the projection morphism $\cU \rw Q$ is flat. 
Now by the universal property of Quot schemes, there is a (universal) quotient sheaf $\cL$ of $\Omega_{\cU/Q}^{**}$ on $\cU$. 
Moreover, by \cite[Th\'{e}or\`{e}me 12.2.1 (v)]{grothendieck1966elements}, the loci 
\[\{q\in Q\, : \,\cL\vert_{\cU_q} \textnormal{ is torsion-free}\} \textnormal{ and } \{q\in Q\, : \,\cL^*\vert_{\cU_q} \textnormal{ is reflexive}\}\]
are open. 
After shrinking $Q$, we may assume that $\cL$ is torsion-free and $\cL^*$ is reflexive over any point of $Q$, $\cL$ is torsion-free, and $\cL^*$ is reflexive. 
In particular, $\Omega_{\cU/Q}^{**}\rw \cL$ gives a family of foliations on $\cU$ over $Q$. 

Now for any $(Y,\sG,D)\in\cS_P$, we dualize the following short exact sequence 
\[\xymatrix{0\ar[r] & \sG \ar[r] & T_Y \ar[r]^-\tau & T_Y/\sG \ar[r] & 0}\] 
and get the short exact sequence 
\[\xymatrix{0\ar[r] & (T_Y/\sG)^* \ar[r]^-{\tau^*} & \Omega_Y^{**} \ar[r] & \textnormal{coker}(\tau^*) \ar[r] & 0.}\] 
Note that 
\[\xymatrix{0\ar[r] & \textnormal{coker}(\tau^*)\otimes\cO_Y(D) \ar[r] & \cO_Y(K_\sG+D) \ar[r] & \cZ\otimes\cO_Y(D) \ar[r] & 0}\] 
where $\cZ$ is the sheaf supported at a $0$-dimensional subscheme. 
Notice that $m(K_\sG+D)+nK_Y$ is ample and Cartier, we have the length of $\cZ$ is at most 
\begin{align*}
& h^0(Y,\cZ\otimes\cO_Y(D+m(K_\sG+D)+nK_Y)) \\
= {} & \chi(Y,\cZ\otimes\cO_Y(D+m(K_\sG+D)+nK_Y)) \\
\leq {} & h^0(Y,\cO_Y((m+1)(K_\sG+D)+nK_Y))\leq\varphi
\end{align*}

Then we have 
\[\chi(Y,\cO_Y((K_\sG+D)+\ell(m(K_\sG+D)+nK_Y))) - \chi(Y,\textnormal{coker}(\tau^*)\otimes\cO_Y(D+\ell(m(K_\sG+D)+nK_Y)))\]
is a non-negative constant $\leq \varphi$. 
Thus, the function 
\[\chi(Y,\textnormal{coker}(\tau^*)\otimes\cO_Y(D+\ell(m(K_\sG+D)+nK_Y)))\] 
has the form $\Phi(\ell)-s$. 
So by the construction of $Q$, there is a $q\in Q$ such that $\rho : Y \rw \cU_q$ is an isomorphism and 
$\textnormal{coker}(\tau^*)\otimes\cO_Y(D)\cong\rho^*(\cL\vert_{\cU_q})$ and thus 
\[\sG(-D) \cong \textnormal{coker}(\tau^*)^*\otimes\cO_Y(-D) \cong (\textnormal{coker}(\tau^*)\otimes\cO_Y(D))^*\cong\rho^*(\cL\vert_{\cU_q})^*.\] 

Notice that $\rho^*(\cL\vert_{\cU_q})^*$ and $\rho^*(\cL^*\vert_{\cU_q})$ are isomorphic at the generic point of $\cU_q$. 
By~\cite[Lemma 1.8]{hacon2021birational}, they are indeed isomorphic. 
Hence, $\cL^*$ gives the desired family of foliations over $Q$. 
\end{enumerate}
\qed
\end{pf}

\begin{rmk}
After some stratification on $Q$, we have $(\cY,\cG,\cD) \rw Q$ is a family of foliated triples (see Defition~\ref{defn_fol_family}) for some reduced divisor $\cD$.
\end{rmk}

\subsection{Boundedness of non-cusp singularities} 
One application of Theorem~\ref{key_thm} is the boundedness of non-cusp singularities. 
Before we show the boundedness of non-cusp singularities, we need the following lemma. 

\begin{lem}\label{nef_big_open}
Let $\cX \rw T$ be a flat family of projective varieties of dimension $n$ and $D$ is a $\bQ$-Cartier divisor on $\cX$ and flat over $T$. 
If there is a $t_0\in T$ such that $D_{t_0}$ is nef (resp. strictly nef) and big, then there exists an open subset $U\subset T$ containing $t_0$ such that $D_t$ is nef (resp. strictly nef) and big for all $t\in U$. 
Moreover, there exists an $N$ such that $h^0(\cX_t,ND_t)\geq 1$ for all $t\in T$. 
\end{lem}
\begin{pf}
We may assume that $T$ is irreducible. 
Indeed, suppose the lemma holds whenever $T$ is irreducible. 
In general, we write $T = \bigcup_{j=1}^mT^j$ as the decomposition of irreducible components $T^j$. 
Since $t_0\in T$, we have $t_0\in T^j$ for some $j$. 
Let us say $t_0\in T^1$. 
By the assumption, we have an open subset $U^1\subset T^1$ containing $t_0$ such that $D_t$ is nef (resp. strictly nef) and big for all $t\in U^1$. 
If $U^1$ is contained in $T\setminus\bigcup_{j=2}^mT^j$ which is an open subset of $T$, then $U^1$ is an open subset of $T$ and we are done. 
If not, there is a point $t_1\in U^1\cap\bigcup_{j=2}^m T^j \subset \bigcup_{j=2}^m T^j$. 
By induction on the number of the irreducible components of $T$, we have an open subset $U^2\subset \bigcup_{j=2}^m T^j$ containing $t_1$ such that $D_t$ is nef (resp. strictly nef) and big for all $t\in U^2$. 
Note that $C^1 := T^1\setminus U^1$ and $C^2 := \bigcup_{j=2}^mT^j\setminus U^2$ are closed subsets of $T$. 
Thus $U := T\setminus(C^1\cup C^2)$ is an open subset of $T$ containing $t_0$ such that $D_t$ is nef (resp. strictly nef) and big for all $t\in U$. 

So from now on, we assume that $T$ is irreducible. 
Since $D_{t_0}$ is nef, we have $D_{t}$ is nef for very general $t\in T^* = T\setminus\bigcup_{i=1}^\infty T_i$ where $T_i$'s are all proper closed subsets of $T$. 
(For a reference, see \cite[Proposition 1.4.14]{lazarsfeld2004positivity}.)
Note that, for all $t\in T^*$, we have $D_t^n = D_{t_0}^n>0$ since $D_{t_0}$ is nef and big.
Hence, $D_t$ is also nef and big for all $t\in T^*$. 

Let $S_N := \{t\in T\vert\, h^0(\cX_t,ND_{t})\geq 1\}$ which is a closed subset in $T$ by the upper semi-continuity. 
Note that we have \[T(\bC) = \left(\bigcup_N S_N(\bC)\right)\cup\left(\bigcup_i T_i(\bC)\right).\] 

By the Baire category theorem, there is an $N$ such that $S_N = T$. 
Thus, the relative base locus of $D$ is proper. 
Let $U$ be an open subset of $T$ such that all irreducible components of the relative base locus of $D$ are flat over $U$. 
Notice that it is possible that $t_0$ is \emph{not} in $U$. 

Now for any $u\in U$, if there is a curve $C_u$ such that $D_u\cdot C_u<0$ (resp. $\leq 0$), then $C_u$ is in the relative base locus of $D$ and moves in the family over $U$. 
So for any $t\in T^*\cap U$, we have $D_t\cdot C_t = D_u\cdot C_u < 0$ (resp. $\leq 0$), which is a contradiction. 
Therefore, $D_u$ is nef (resp. strictly nef) for all $u\in U$.  
Moreover, $D_u^n = D_{t}^n >0$ where $t\in T^*\cap U$, then $D_u$ is also big for $u\in U$. 

If $t_0\in U$, then we are done. 
If not, that is $t_0\in T\setminus U$, then by induction on $\dim T$, there exists an open subset $V\subset T\setminus U$ containing $t_0$ such that $D_t$ is nef (resp. strictly nef) and big for all $t\in V$. 
So $(T\setminus U)\setminus V$ is a closed subset of $T\setminus U$, which is a closed subset of $T$. 
Thus, $(T\setminus U)\setminus V$ is a closed subset of $T$. 
Therefore $V\cup U = T\setminus\big((T\setminus U)\setminus V\big) \subset T$ is an open subset containing $t_0$ such that $D_t$ nef (resp. strictly nef) and big for all $t\in V\cup U$.
\qed
\end{pf}

This proof also gives the following Corollary. 
\begin{cor}\label{trivial_move}
Every $(K_\sG+D)$-trivial curve moves in the family over $U$. 
\end{cor}

\begin{prop}\label{dihedral_bdd}
Fix a function $P:\bZ_{\geq 0} \rw \bZ$. 
There exist constants $C_1$, $C_2$, and $C_3$ such that, for any log canonical model $(X,\sF)$ with Hilbert function $\chi(X,mK_{\sF})=P(m)$, the index (resp. embedding dimension) at the dihedral singularities is at most $C_1$ (resp. $C_2$) and the number of dihedral singularities is at most $C_3$. 
Moreover, there exist constants $B_1$, $B_2$, and $B_3$ such that the index (resp. embedding dimension) at the non-terminal foliation singularities which are cyclic quotient singularities is at most $B_1$ (resp. $B_2$) and the number of non-terminal foliation singularities which are cyclic quotient singularities is at most $B_3$. 
\end{prop}
\begin{pf}
The proof follows closely from the proof of \cite[Proposition 3.1]{chen2021boundedness} for canonical models. 
For the reader's convenience, we include the proof. 

By Theorem~\ref{key_thm}, there exists a projective morphism $\mu : \cY\rw T$ of quasi-projective varieties of finite type and a family of foliations $\cF$ on $\cY$ over $T$ such that for any $(Y,\sG,D)\in\cS_P$, there exists a $t\in T$ and an isomorphism $\rho:Y\rw \cY_t$ such that $\sG(-D)\cong \rho^*(\cF\vert_{\cY_t})$. 
We may assume that 
\[T^* :=\{t\in T\,\vert\, (\cY_t,\cF_t)\cong (Y,\sG(-D)) \mbox{ for some } (Y,\sG,D)\in\cS_P\}\] 
is dense in $T$. 
Fix a $t_0\in T^*$. 

Note that $K_\cF\vert_{\cY_t}$ is nef and big, by Lemma~\ref{nef_big_open}, we have that $K_\cF$ is big over $T$. 
Thus, the relative base locus of $K_\cF$ is proper. 
Therefore, we have a finite stratification of $T$ such that all irreducible components of the relative base locus of $K_\cF$ are flat over $T$. 

Let $\cE_j$'s be those irreducible components of the relative base locus of $K_\cF$ over $U$ of codimension one. 
Note that, for any $t\in T^*$ and for any irreducible curve $E\subset \cY_t$ with $K_\cF\cdot E = 0$, we have $E$ is in the relative base locus of $D$ and hence $E\subset\cE_j$ for some $j$. 
Let $\pi:\wt{\cY}\rw\cY$ be the log resolution of $(\cY,\sum_j\cE_j)$ with the exceptional divisors $\sum \cF_\ell$. 
By shrinking $T$ further, we may assume that $\nu := \mu\circ \pi$ is smooth and $\cF_\ell$ is flat over $T$. 
Let $\wt{\cE_j}$ be the proper transform of $\cE_j$ via $\pi$. 

Note that $\Lambda := \{\wt{\cE_j}, \cF_\ell\}$ is a finite set. 
So there are only finitely many possible connected subsets of $\Lambda$ with dual graphs of $D_n$ type or $A_n$ type. 
By \cite[Satz 2.9 and Satz 2.11]{brieskorn1968rationale}, the intersection matrix of a dihedral singularity determines the embedding dimension and the group acting on it. 
Since the intersection matrix remains the same in the family, we have the boundedness of the embedding dimension and the order of the group acting on these singularities. 
Hence the index of $X$ at the dihedral singularities is also bounded. 
Similarly, by \cite{reid2012surface}, the embedding dimension and the index of a cyclic quotient singularity is determined by the intersection matrix. 
Therefore, we have the boundedness of the embedding dimension and of the index at the cyclic quotient singularities. 
\qed
\end{pf}

\begin{rmk}
\begin{enumerate}
\item It is possible that some connected subset of $\Lambda$ of $D_n$ type or of $A_n$ type is redundant. 
\item By the same argument as above, we also show the boundedness of the cusps. 
Precisely, the number of exceptional divisors and the intersection matrix over any cusp singularity are bounded. 
\end{enumerate}
\end{rmk}

We also have the following theorem. 
\begin{thm}\label{bir}
Fix a function $P: \bZ_{\geq 0} \rw \bZ$. 
There exists an integer $m_P$ such that if $(X,\sF)$ is a log canonical model of general type and $\chi(mK_\sF)=P(m)$ for all $m\geq 0$, then for any $m>0$ divisible by $m_P$, $|mK_\sF|$ defines a birational map which is an isomorphism on the complement of the cusp singularities. 
\end{thm}
\begin{pf}
This follows closely from \cite[Theorem 4.3]{hacon2021birational}. 
See also \cite[Theorem 3.3]{chen2021boundedness}. 
For the reader's convenience, we include the proof. 

Let $g : (Y,\sG) \rw (X,\sF)$ be the minimal resolution of cusps. 
Then $K_\sG = g^*K_\sF$ and $i(\sG)=i_\bQ(\sF)$. (See Definition~\ref{index}.)
\begin{claim}
$K_Y+3i(\sG)K_\sG$ is nef. 
\end{claim}
\begin{pf}[Claim]
Given any irreducible curve $C$ on $Y$. 
\begin{enumerate}
\item If $C$ is contracted by $g$, then 
\[(K_Y+3i(\sG)K_\sG)\cdot C = K_Y\cdot C = -2-C^2 \geq 0\]
since $C^2\leq -2$ by the minimality of $g$. 
\item If $C$ is not contracted by $g$, 
\begin{enumerate}
\item If $K_Y\cdot C\geq 0$, then 
\[(K_Y+3i(\sG)K_\sG)\cdot C \geq 3i(\sG)K_{\sF}\cdot g_*C > 0\] 
by Lemma~\ref{ample}.

\item If $K_Y\cdot C<0$, by \cite[Proposition 3.8]{fujino2012minimal}, every $K_Y$-negative extremal ray is spanned by a rational curve $C$ with $0<-K_Y\cdot C\leq 3$. 
Thus, if we write $C \equiv \sum a_iC_i+D$ where $D\cdot K_Y\geq 0$, the $C_i$'s are $K_Y$-negative extremal rays, the $a_i$'s are non-negative, and at least one of $a_i$ is positive. 
By what we have seen above, $(K_Y+3i(\sG)K_\sG)\cdot D\geq 0$. 
Thus, 
\begin{align*}
& (K_Y+3i(\sG)K_\sG)\cdot C \\
= {} & \sum_ia_i \Big((K_Y+3i(\sG)K_\sG)\cdot C_i\Big) + (K_Y+3i(\sG)K_\sG)\cdot D \\
\geq {} & \sum_ia_i(- 3+3) + 0 = 0.
\end{align*}
\end{enumerate} 
\end{enumerate}
This completes the proof of the nefness of $K_Y+3i(\sG)K_\sG$. 
\end{pf}

Now, by \cite[Lemma 3.5]{hacon2021birational}, we know that $\gamma K_\sG-K_Y$ is big for any fixed integer $\gamma$ with 
\[\gamma > \left\lceil\max\left\{\frac{2K_{\sF}\cdot K_{X}}{K_{\sF}^2}+3i_\bQ(\sF),0\right\}\right\rceil.\]

Note that $L := (\beta-\gamma)K_\sG+(\gamma K_\sG-K_Y)$ is pseudo-effective for any $\beta>\gamma$.
Let $L=P+N$ be the Zariski decomposition. 
Then we have 
\[P^2\geq(\beta-\gamma)^2K_\sG^2\geq \left(\frac{\beta-\gamma}{i_\bQ(\sF)}\right)^2.\]
And if $C$ is not in the exceptional divisor of $g$, then 
\[P\cdot C\geq (\beta-\gamma)K_\sG\cdot C\geq \frac{\beta-\gamma}{i_\bQ(\sF)}.\]

By Proposition~\ref{dihedral_bdd} and Proposition~\ref{bdd}, we have the boundedness of embedding dimension of each rational singularity on $Y$. 
Thus, we have a uniform bound $\delta$ such that $\delta_\zeta\leq\delta$ for $|\zeta|=2$. 
Now we fix a $\beta$ divisible by $i_\bQ(\sF)$ such that $P^2>\delta$ and $P\cdot C>\delta/2$ for any irreducible curve $C$. 
By \cite[Theorem 0.1]{langer2001adjoint} and Proposition~\ref{sep}, we have that $K_Y+L = \beta K_\sG$ is very ample on the complement of the exceptional divisor of $g$. 
By \cite[Theorem 6.2]{sakai1984weil}, $g_*\cO_Y(\beta K_\sG) = \cO_X(\beta K_\sF)$. 
Hence, $\beta K_\sF$ is very ample on the complement of the cusps.
\qed
\end{pf}

\subsection{{The moduli functor \texorpdfstring{$\cM_P^{\tn{sm}}$}{M\_Psm}}}\label{sm_stable}
We will define a moduli functor parametrizing stable $\cS_P$-smoothable foliated triples. 
\begin{defn}\label{defn_fol_family}
Fix $P : \bZ_{\geq 0}\rw \bZ$. 
We call $\pi : (\cX,\cF,\Delta) \rw B$ \emph{a family of foliated triples} (see also Definition~\ref{fol_triple}) if 
\begin{enumerate}
\item $\pi : \cX\rw B$ is a projective morphism such that all fibers are surfaces with at worst semi-log canonical singularities. (For a reference, see~\cite[Definition 3.13.5]{hacon2010classification}.)
\item $\cF$ is a foliation of rank 1 on $\cX$ such that $\cF$ is tangent to the fibers of $\pi$, that is, $\cF$ is a family of foliations on $\cX$ over $B$. (See Definition~\ref{family_fol}.)
\item $\Delta = \sum\Delta_i$ is a reduced divisor with all $\Delta_i$'s flat over $B$. 
\end{enumerate}
\end{defn}

\begin{defn}\label{defn_moduli}
Fix an integral-valued function $P : \bZ_{\geq 0} \rw \bZ$. 
\begin{enumerate}
\item A triple $(X,\sF,D)$ is called \emph{$\cS_P$-smoothable} if the followings hold:
\begin{enumerate}
\item $X = \cup_iX_i$ is a surface with at worst semi-log canonical singularities where $X_i$'s are irreducible. 
\item $\sF$ is a pre-foliation on $X$. 
\item $D$ is a Weil divisor on $X$. 
\item There exists a flat projective $\bQ$-Gorenstein (that is, both $K_\cX$ and $K_\cF$ are $\bQ$-Cartier) one-parameter family of foliated triples $\pi : (\cX,\cF,\D) \rw C$ such that, for some point  $0\in C$, the restriction on the fiber over the point $c\in C\setminus\{0\}$ is in $\cS_P$ and the restriciton on the fiber $\cX_0$ over $0$ is $(X,\sF,D)$, that is, $X\cong\cX_0$, $\cF\vert_{\cX_0}\cong\sF$, and $\D\vert_{\cX_0}=D$. 
\item $(\cF,\Delta)$ and $(\cF_{X_i},\Delta_{X_i})$ are log canonical for all $i$ where $\cF_{X_i}$ is the restricted foliation and $K_{\cF_{X_i}}+\Delta_{X_i} = (K_\cF+\Delta)\vert_{X_i}$. 
\end{enumerate}

\item An $\cS_P$-smoothable triple $(X,\sF,D)$ is called \emph{stable} if $\Omega_\cF(\Delta)\vert_X$ is nef and 
\[\left(\left(\omega_\cX\otimes(\Omega_\cF(\Delta))^{4i(P)}\right)^{**}\right)\vert_X\] 
is ample where $\pi : (\cX,\cF,\Delta)\rw C$ is a flat projective $\bQ$-Gorenstein one-parameter family of foliated triples as in  (1)(d).
(See Lemma~\ref{iP} for the definition of $i(P)$.)

\item A family of \emph{stable $\cS_P$-smoothable} foliated triples over a scheme $B$ is a family of foliated triples $\pi : (\cX,\cF,\Delta) \rw B$ over $B$ such that, for any $b\in B$, 
\begin{enumerate}
\item the fiber $(Y,\cF,\Delta)_b$ is a \emph{stable $\cS_P$-smoothable} foliated triple, 
\item $\omega_{\cX/B}^{[\ell]}\vert_{\cX_b}\cong \omega_{\cX_b}^{[\ell]}$ for all $\ell\in\bZ$, and
\item $(\Omega_{\cF}(\Delta))^{[\ell]}\vert_{\cX_b}\cong (\Omega_{\cF}(\Delta)\vert_{\cX_b})^{[\ell]}$ for all $\ell\in\bZ$ where $(\bullet)^{[\ell]} = ((\bullet)^{\otimes \ell})^{**}$. 
\end{enumerate}
\end{enumerate}
\end{defn}

\begin{rmk}
In particular, any foliated triples $(Y,\sG,D)$ in $\cS_P$ are stable $\cS_P$-smoothable. 
Indeed, we consider $(X,\sF)$ the log canonical model of $(Y,\sG,D)$ (see Remark~\ref{unique}) and $f : (Y,\sG,D) \rw (X,\sF)$ the associated morphism. 
In the proof of Theorem~\ref{key_thm}, we have shown that $K_\sG+D = f^*K_\sF$ which is nef and $K_Y+mi(\sG)(K_\sG+D)$ is ample for $m\geq 4$. 
\end{rmk}

\begin{rmk}
Let $X$ be a surface of a threefold. 
Given any log canonical foliated triple $(X,\sF,D)$ of rank one. 
For any proper birational morphism $\pi : Y \rw X$, we have 
\[K_{\sF_Y}+\pi_*^{-1}D+F = \pi^*(K_\sF+D)+E\]
where $E$ and $F$ are effective exceptional divisors and $F$ is reduced by \cite[Lemma 8.2]{cascini2020mmp}. 
The \emph{pullback foliated triple} is defined to be $(Y,\sF_Y,\pi_*^{-1}D+F)$. 
In particular, we have $\pi_*(K_{\sF_Y}+\pi_*^{-1}D+F) = K_\sF+D$. 
\end{rmk}

\begin{defn}\label{functor}
We define the moduli functor $\cM_P^{\tn{sm}} :  \tn{(Schemes)} \rw \tn{(Sets)}$ by sending the scheme $B$ to the sets of isomorphism classes of \emph{stable $\cS_P$-smoothable} foliated triples over $B$. 
\end{defn}

\section{Invariance of plurigenera}
In this section, we would like to show the invariance of plurigenera for the family of minimal partial du Val resolutions of log canonical models of general type with fixed Hilbert function. 

\begin{lem}\label{iP}
Fix a function $P : \bZ_{\geq 0} \rw \bZ$. 
Then there exists an $i(P)\in\bN$ depending only on $P$ such that $i(P)(K_\sG+D)$ is Cartier for all $(Y,\sG,D)\in\cS_P$. 
\end{lem}
\begin{pf}
We can take $i(P) = C_2!$ where $C_2$ comes from Proposition~\ref{bdd} and depends only on $P$. \qed
\end{pf}

\begin{thm}\label{inv}
Fix a function $P : \bZ_{\geq 0} \rw \bZ$. 
Given any family of foliated triples $f : (\cX,\cF,\Delta)\rw T$ where $T$ is connected and $(\cX,\cF,\Delta)_t\in\cS_P$ for all $t\in T$ (see Definition~\ref{defn_fol_family} for families of foliated triples), there exists an $m_{f}$ depending only on this family $f$ such that  
\[h^0(\cX_t,m(K_{\cF}+\Delta)\vert_{\cX_t})\] 
is independent of $t$ for $m\geq m_f$ and divisible by $i(P)$. (See Lemma~\ref{iP}.)
\end{thm}

\begin{pf}
The idea of the proof is from the proof of \cite[Proposition 3.7]{cascini2018invariance}. 
We divide the proof into the following steps. 
Let $m\in\bN$ be a positive integer divisible by $i(P)$. 
\begin{enumerate}
\item Fix a point $0\in T$ and let $(X,\sF,D) := (\cX,\cF,\Delta)\vert_{\cX_0}\in\cS_P$. 
Let $\rho : X^{\tn{m}} \rw X$ be the minimal resolution of $X$ and $\sF^{\tn{m}}$ be the pullback foliation on $X^\tn{m}$. 
Contracting down all $\sF^{\tn{m}}$-chains to get $(Y,\sG)$, we have a morphism $\pi: (Y,\sG) \rw (X,\sF)$. 
\[\xymatrix{X^{\tn{m}}  \ar[dd]_-\rho \ar[rd] & \\ & Y \ar[ld]^-\pi \\ X &}\]
Let $E$ be the proper transform of $D$. 
Thus, because of the construction of the minimal partial du Val resolution (\ref{mpdvr}), we have $K_\sG+E = \pi^*(K_\sF+D)$, which is nef and big.  

Let $\Gamma_1$, $\ldots$, $\Gamma_p$ be all elliptic Gorenstein leaves on $Y$ and $C_1$, $\ldots$, $C_q$ be all disconnected chains of rational curves \emph{exceptional over its log canonical model} such that $(K_\sG+E)\cdot C_j^{(k)}=0$ where $C_j = \sum_{k=1}^{r_j} C_j^{(k)}$ is the sum of its irreducible components. 
Also let $Z = \sum_{i=1}^p\Gamma_i + \sum_{j=1}^q C_j$. 

\item \textit{Claim.} $h^1(\Gamma_i,m(K_{\sG}+E))=0$ for all $i$.
\begin{pf}
Notice that $\Gamma_i$ is Cohen-Macaulay with the trivial dualizing sheaf. 
Since each irreducible component of $E$ is over a log canonical non-canonical foliation singularity, $E$ does not meet $\Gamma_i$. 
So we have 
\[h^1(\Gamma_i, m(K_{\sG}+E)) = h^0(\Gamma_i,-m(K_{\sG}+E)) = h^0(\Gamma_i,-mK_{\sG}) = 0\]
where the last equality is because $\cO_{\Gamma_i}(mK_{\sG})$ is not torsion and has degree 0 by \cite[FactIII.0.4 and Theorem IV.2.2]{mcquillan2008canonical}.
\end{pf}

\item \textit{Claim.} $h^1(C_j,m(K_{\sG}+E))=0$ for all $j$.
\begin{pf}
Note that we have  
\[(K_\sG+E)\cdot C_j^{(k)}=0\]
for all $j$ and $k$. 
Since $C_j^{(k)}$'s are either smooth rational curves or a rational curve with two singularities of index two, we have 
\begin{equation}\label{h1_van}
h^1(C_j^{(k)},m(K_\sG+E)) = 0
\end{equation}
for all $j$ and $k$. 
Now, for any fixed $j$, we consider the following short exact sequence 
\[\xymatrix{0 \ar[r] & \cO_{C_j} \ar[r] & \displaystyle \bigoplus_{k=1}^{r_j}\cO_{C_j^{(k)}} \ar[r] & \displaystyle \bigoplus_{k=1}^{r_j-1}\cO_{C_j^{(k)}\cap C_j^{(k+1)}} \ar[r] & 0.}\]
Tensoring with $\cO_{C_j}(m(K_\sG+E))$ and taking the induced long exact sequence, we notice that the map 
\[\xymatrix{\displaystyle \bigoplus_{k=1}^{r_j}H^0\big({C_j^{(k)}},m(K_\sG+E)\big) \ar[r] & \displaystyle \bigoplus_{k=1}^{r_j-1}H^0\big({C_j^{(k)}\cap C_j^{(k+1)}},m(K_\sG+E)\big)}\]
is surjective and thus, we have that $h^1(C_j,m(K_{\sG}+E))=0$ for all $j$ by equation~(\ref{h1_van}). 
\end{pf}

\item \textit{Claim.} $h^j(Y, m(K_{\sG}+E)-Z)=0$ for $j>0$ and for all $m\geq m_0$ for some $m_0$. 
\begin{pf}
For any irreducible curve $C$ in the support of $Z$, we have that 
\[\big(m(K_{\sG}+E)-(K_Y+Z)\big)\cdot C = -(K_Y+Z)\cdot C\geq 0.\]
Note that $K_\sG+E = \pi^*(K_\sF+D)$ is nef and big, and the $(K_\sG+E)$-trivial curves are exactly the curves in the support of $Z$.
Then there exists an $m_0$ depending only on $(X,\sF,D)$ such that $m(K_{\sG}+E)-(K_Y+Z)$ is nef and big for $m\geq m_0$. 
Since $Y$ is klt, we have 
\[h^j(Y,m(K_{\sG}+E)-Z) = h^j(Y,K_Y+m(K_{\sG}+E)-(K_Y+Z)) = 0\]
for all $j>0$ and $m\geq m_0$ by Kawamata-Viehweg vanishing theorem. 
\end{pf}

\item \textit{Claim.} $h^j(Y,m(K_{\sG}+E))=0$ for $j>0$ and for all $m\geq m_0$. 
\begin{pf}
By steps (2) and (3), we have that $h^1(Z,m(K_\sG+E)) = 0$. 
Now we consider the short exact sequence 
\[\xymatrix{0 \ar[r] & \cO_Y(-Z) \ar[r] & \cO_Y \ar[r] & \cO_Z \ar[r] & 0.}\]
Tensoring with $\cO_Y(m(K_\sG+E))$ and taking the induced long exact sequence, we have the map
\[\xymatrix{H^1(Y,m(K_{\sG}+E)-Z) \ar[r] & H^1(Y,m(K_{\sG}+E))}\]
is surjective and $H^2(Y,m(K_{\sG}+E)-Z) \cong H^2(Y,m(K_{\sG}+E))$. 
Thus, we have $h^j(Y,m(K_{\sG}+E)) = 0$ for $j>0$ and for all $m\geq m_0$ by step (4). 
\end{pf}

\item \textit{Claim.} $h^1(X,m(K_\sF+D))=0$ for all $m\geq m_0$. 
\begin{pf}
Note that $\pi_*(m(K_\sG+E)) = m(K_\sF+D)$ by the projection formula. 
Since $h^1(Y,m(K_{\sG}+E))=0$ for all $m\geq m_0$ by step (5), we have that 
\[h^1(X,m(K_{\sF}+D))=0\] 
for all $m\geq m_0$ by Leray spectral sequence. 
\end{pf}

\item \textit{Claim.} $h^2(\cX_t,m(K_{\cF_t}+\Delta_t))=0$ for all $t\in T$ and for $m\geq m_P$ for some $m_P$. 
\begin{pf}
Note that $K_{\cF_t}+\Delta_t$ is nef and big since $(\cX_t,\cF_t,\Delta_t)\in\cS_P$.  
In addition, we have already seen, in the proof of Theorem~\ref{key_thm}, that $K_{\cX_t}+4i(P)(K_{\cF_t}+\Delta_t)$ is ample. 
Moreover, by \cite[Lemma 3.5]{hacon2021birational}, we have that $m(K_{\cF_t}+\Delta_t) - K_{\cX_t}$ is big for $m\geq m_P$ where 
\[m_P := \max\left\{\frac{2(K_{\cF_t}+\Delta_t)\cdot K_{\cX_t}}{(K_{\cF_t}+\Delta_t)^2}+4i(P)+1, 1\right\}.\]
By Proposition~\ref{bdd}, we have that $m_P$ depends only on $P$. 

Now by Serre duality, we have, for $m\geq m_P$, that 
\[h^2(\cX_t,m(K_{\cF_t}+\Delta_t)) = h^0(\cX_t,K_{\cX_t} - m(K_{\cF_t}+\Delta_t)) = 0\] 
since there is no section for an anti-big divisor. 
\end{pf}

\item Take $m_f := \max\{m_0, m_P\}$ which only depends on the family $f$. 
Now we have, for all $t\in T$ and for all $m\geq m_f$, that  
\begin{align*}
h^0(\cX_t,m(K_{\cF_t}+\Delta_t)) &\geq \chi(\cX_t,m(K_{\cF_t}+\Delta_t)) \\
&= \chi(\cX_0,m(K_{\cF_0}+\Delta_0)) \\
&= h^0(\cX_0,m(K_{\cF_0}+\Delta_0))
\end{align*}
where the first inequality comes from step (7) and the last equality comes from steps (6) and (7). 

On the other hand, by upper semi-continuity, we have 
\[h^0(\cX_t,m(K_{\cF_t}+\Delta_t))\leq h^0(\cX_0,m(K_{\cF_0}+\Delta_0)).\]
Therefore, we get the desired equality. 
\qed
\end{enumerate}
\end{pf}

\begin{rmk}\label{h_vanishing}
The arguments above also show that $h^1(\cX_t,m(K_{\cF_t}+\Delta_t)) = 0$ for $m\geq m_f$ and divisible by $i(P)$. 
\end{rmk}

\begin{cor}
Fix a function $P : \bZ_{\geq 0} \rw \bZ$. 
There is an $m_P$ depending only on $P$ such that, for all $m\geq m_P$ and divisible by $i(P)$, we have that 
\[h^1(X,m(K_\sF+D))=h^2(X,m(K_\sF+D))=0\] and 
$h^0(X,m(K_\sF+D))$ is constant for all $(X,\sF,D)\in\cS_P$. 
\end{cor}
\begin{pf}
By Theorem~\ref{key_thm} and \ref{open}, we have finitely many families $f_j$'s of foliated triples parametrizing $\cS_P$. 
Taking $m_P$ divisible by all $m_{f_j}$'s in Theorem~\ref{inv}, we have 
\[h^1(X,m(K_\sF+D)) = h^2(X,m(K_\sF+D)) = 0\] 
for all $(X,\sF,D)\in\cS_P$. 
Since $\chi(X,m(K_\sF+D))$ is constant for all $(X,\sF,D)\in\cS_P$, so is $h^0(X,m(K_\sF+D))$. 
\qed
\end{pf}

\section{{Valuative criterion of properness for \texorpdfstring{$\cM_P^{\tn{sm}}$}{M\_Psm}}}
In this section, we establish the relative $(K_\cF+\D)$-minimal model program for the families of foliated triples which is a semi-stable reduction of a family of foliated triples whose general fiber is in $\cS_P$. 
Also we show the valuative criterion of properness for the moduli functor $\cM_P^{\tn{sm}}$ for families of foliated triples whose general fiber is in $\cS_P$ (see Theorem~\ref{properness}).

\subsection{Elementary singularities}
\begin{defn}\label{elementary}
Let $X$ be a smooth variety, $\sF$ be a foliation of rank 1 on $X$, and $p$ be a singularity of $\sF$. 
Let $v$ be a germ of a vector field around $p$ associated to $\sF$. 
Then $p$ is called \emph{elementary} if $v$ has a non-zero non-nilpotent linear part. 
\end{defn}

\begin{rmk}
Elementary singularities are log canonical by \cite[Fact I.ii.4]{mcquillan2013almost}. 
\end{rmk}

By the Introduction and Theorem 1 in \cite{cano2014vector}, we have the following. 
\begin{thm}[{\cite[Introduction and Theorem 1]{cano2014vector}}]\label{line_rsln}
Any germ of vector fields on $(\bC^3,0)$ tangent to a rank-two foliation can be desingularized. 
That is, there exists a finite sequence of blowups such that the pullback foliation has only elementary singularities. 
\end{thm}

\begin{defn}[{\cite[Definition 7.1]{kollar1998birational}}]\label{dlt}
Let $X$ be a normal variety, $B$ an effective $\bQ$-divisor, and $f: X\rw C$ a non-constant morphism to a smooth curve $C$. 
We say that $f : (X,B) \rw C$ is a \emph{divisorial log terminal} (dlt) if the pair $(X,B+f^{-1}(c))$ is dlt for every closed point $c\in C$. 
\end{defn}

\subsection{Semi-stable reduction}\label{ss_reduction}
Given a family of foliated triples $f : (\overline{\cX},\overline{\cF},\overline{\Delta}) \rw C$ over a smooth curve $C$ whose general fiber is in $\cS_P$. 
After taking a base change $T \rw C$, we have a semi-stable reduction for the family of surfaces $(\overline{\cX},\overline{\Delta})$ over the curve $C$ by \cite{kempf1973toroidal} (see also \cite[Theorem 7.17]{kollar1998birational}) and we denote this family as $f : (\cX,\cF,\Delta)\rw T$. 

Since $(\cX,\Delta)$ is log smooth and $\cF$ is tangent to the fibration which is a rank-two foliation, by Theorem~\ref{line_rsln}, we may assume, after a finite sequence of blowups, that $\cF$ has only elementary singularities and each log canonical center of $\cF$ is not contained in $\Delta$. 

Moreover, by taking base changes and resolutions, we may assume the central fiber of $f$ is reduced. 
We call $f : (\cX,\cF,\Delta) \rw T$ a semi-stable reduction of the family of foliated triples $\overline{f}:(\overline{\cX},\overline{\cF},\overline{\Delta})\rw C$ over $C$. 
By the construction, $(\cX,\cF,\D)$ is foliated log smooth. 
In particular, $f : \cX \rw T$ is dlt (see Definition~\ref{dlt}). 

\begin{thm}\label{KF_MMP}
Let $f : (\cX,\cF,\D) \rw T$ be a family of foliated triples which is a semi-stable reduction of a family of foliated triples whose general fiber is in $\cS_P$. 
Then the relative $(K_\cF+\D)$-minimal model program exists and stops with a family of foliated triples $(\cX^{\tn{m}},\cF^{\tn{m}},\Delta^{\tn{m}}) \rw T$ such that $(K_{\cF^{\tn{m}}}+\D^{\tn{m}})$ is nef over $T$, $f^{\tn{m}} : \cX^{\tn{m}} \rw T$ is dlt, and the general fiber is in $\cS_P$. 
\end{thm}
\begin{pf}
Let $R$ be a $(K_\cF+\Delta)$-negative extremal ray over $T$. 
Then its locus is not the whole space $\cX$ since the general fiber of $f$ is a model of a foliated triple in $\cS_P$. 
Thus, by the proof of \cite[Theorem 8.8]{cascini2020mmp}, we may run the $(K_\cF+\Delta)$-MMP over $T$, which will end up with a family of foliated triples $(\cX^{\tn{m}},\cF^{\tn{m}},\Delta^{\tn{m}}) \rw T$ such that $K_{\cF^{\tn{m}}}+\D^{\tn{m}}$ is nef over $T$. 
Moreover, $f^{\tn{m}} : \cX^{\tn{m}} \rw T$ is dlt and the general fiber is in $\cS_P$. 
\qed
\end{pf}

Next we show the valuative criterion of properness for the moduli functor $\cM_P^{\tn{sm}}$ for families whose general fiber is in $\cS_P$. 

\begin{thm}\label{properness}
Let $T$ be a smooth curve and $(\cX^\circ,\cF^\circ,\D^\circ)\in\cM_P^{\tn{sm}}(T^\circ)$ be a family of stable $\cS_P$-smoothable foliated triples over $T^\circ$ where $T^\circ = T\setminus\{0\}$ and $0\in T$. 
Suppose the fiber of the family over $t\neq 0$ is in $\cS_P$. 
Then there exists a smooth curve $C$ and a finite morphism $C \rw T$ such that, after the base change via $C \rw T$, the pullback family $(\cX_C^\circ,\cF_C^\circ,\D_C^\circ)$ extends uniquely to a family of stable $\cS_P$-smoothable foliated triples $(\cX_C,\cF_C,\D_C)$ over $C$. 
\end{thm}
\begin{pf}
Let $g : (\cY,\cG,\T) \rw C$ be a semi-stable reduction of the closure of the family $f^\circ : (\cX^\circ,\cF^\circ,\D^\circ) \rw T^\circ$ where $C \rw T$ is a finite morphism and $C$ is a smooth curve. 

Then by Theorem~\ref{KF_MMP}, we run the relative $(K_\cG+\T)$-MMP over $C$. 
Thus, we get a family $g : (\cY,\cG,\T) \rw C$ such that $K_\cG+\T$ is nef over $C$. 
By running additional partial $K_\cY$-MMP over $C$, we may assume that $K_\cY+4i(P)(K_\cG+\T)$ is nef over $C$. 
(See Lemma~\ref{iP} for the definition of $i(P)$.)
Note that $\big(K_\cY+4i(P)(K_\cG+\T)\big)-K_\cY$ is nef and big over $C$. 
By relative base point free theorem (see \cite[Theorem 3.24]{kollar1998birational} for a reference), $K_\cY+4i(P)(K_\cG+\T)$ is semi-ample over $C$. 
Thus, we have the relative canonical model $g : (\cY,\cG,\T) \rw C$. 

Let $f^\circ_C$ be the base change of $f^\circ$ via the morphism $C^\circ \rw T^\circ$ where $C^\circ = C\setminus\{0\}$. 
Notice that for any $c\in C\setminus\{0\}$, the fibers of $f^\circ_C$ and $g$ over $c$ are isomorphic. 
Moreover, for $c=0$, we have that $(K_\cG+\T)\vert_{\cY_0}$ is nef and $K_{\cY_0}+4i(P)(K_{\cG}+\T)\vert_{\cY_0}$ is ample. 
Hence, the fiber of $g$ over $0\in C$ is stable and $\cS_P$-smoothable. 
Therefore, $g : (\cY,\cG,\T) \rw C$ is a family of stable $\cS_P$-smoothable foliated triples over $C$. 
\qed
\end{pf}

\section{{Valuative criterion of separatedness for \texorpdfstring{$\cM_P^{\tn{sm}}$}{M\_Psm}}}
In this section, we show the following theorem. 

\begin{thm}\label{separatedness}
Let $T$ be a smooth curve, $T^\circ = T\setminus\{0\}$ where $0\in T$, and $(\cX_i,\cF_i,\D_i) \in\cM_P^{\tn{sm}}(T)$ be two families of stable $\cS_P$-smoothable foliated triples over $T$ for $i=1$, $2$. 
Then any isomorphism $\varphi^\circ : (\cX_1,\cF_1,\D_1)_{T^\circ} \rw (\cX_2,\cF_2,\D_2)_{T^\circ}$ over $T^\circ$ extends uniquely to an isomorphism $\varphi : (\cX_1,\cF_1,\D_1) \rw (\cX_2,\cF_2,\D_2)$ over $T$. 
\end{thm}

The proof follows the classical approach. 
We will need the following lemmas that, for $(Y,\sG,D)\in\cS_P$, $Y$ is indeed a Proj of a (finitely generated) graded ring associated with the adjoint bundle $K_Y+4i(P)(K_\sG+D)$ and this adjoint bundle behaves well after taking resolutions.

\begin{lem}
Fix a function $P : \bZ_{\geq 0} \rw \bZ$. 
For any $(Y,\sG,D)\in\cS_P$, we have 
\[Y \cong \textnormal{Proj}\bigoplus_{m\geq 0} H^0(Y,m(K_Y+4i(P)(K_\sG+D)))\] 
where $i(P)$ is defined in Lemma~\ref{iP}.
\end{lem}
\begin{pf}
Notice that $K_Y+4i(P)(K_\sG+D)$ is ample by the part (1) of the proof of Theorem~\ref{key_thm}. \qed
\end{pf}

\begin{lem}\label{sep_lc}
Fix a function $P : \bZ_{\geq 0} \rw \bZ$. 
For any $(Y,\sG,D)\in\cS_P$ and any resolution of singularities $f : Z \rw Y$, we write 
\[K_Z+4i(P)\left(K_\sH+f_*^{-1}D+\sum_j\ve(E_j)E_j\right) = f^*(K_Y+4i(P)(K_\sG+D))+\sum_j a_jE_j\] 
where $\sH$ is the pullback foliation on $Z$ and $E_j$'s are exceptional divisors of $f$. 
Then we have $a_j > -1$ for all $j$ where $\ve(E_j) = 0$ if $E_j$ is $\sH$-invariant and $\ve(E_j)=1$ otherwise.  
\end{lem}
\begin{pf}
Let $(X,\sF)$ be the log canonical model of $(Y,\sG,D)$. 
Note that $K_\sG+D$ is the pullback of $K_\sF$. 
Thus the divisor $K_\sH+f_*^{-1}D+\sum_j\ve(E_j)E_j - f^*(K_\sG+D) = \sum_ja'_jE_j$ is effective and exceptional since $\sF$ has at worst log canonical singularities. 
Moreover, since $Y$ has at worst klt singularities, we have $K_Z - f^*K_Y = \sum_j a''_jE_j$ where $a''_j > -1$. 
Hence, $a_j = 4i(P)a'_j+a''_j > -1$. 
\qed
\end{pf}

\begin{lem}\label{adjunction}
Fix a function $P : \bZ_{\geq 0} \rw \bZ$. 
Let $T$ be a smooth affine curve and $(\cX,\cF,\D)\in\cM_P^{\tn{sm}}(T)$ be a family of stable $\cS_P$-smoothable foliated triples over $T$. 
Suppose the fiber of the family over $t\neq 0$ is in $\cS_P$. 
Assume there is a birational morphism $g : \cY \rw \cX$ such that $f\circ g : \cY \rw T$ is 
a semi-stable reduction of $(\cX,\cF,\D)$. 
Let $\cG$ be the pullback foliation of $\cF$ and $\T$ be the proper transform of $\D$. 
We may write
\[K_\cY+4i(P)(K_\cG+\T) = g^*\big(K_\cX+4i(P)(K_\cF+\D)\big) + \sum_i a_iE_i + \sum_j b_jF_j\]
where $E_i$ are all exceptional divisors mapping to $0\in T$ and $F_j$ are all exceptional divisors flat over $T$. 
Then $a_i\geq 0$ for all $i$ and $b_j\geq -1$ for all $j$. 
\end{lem}
\begin{pf}
Restriction the equality to the general element $t\in T$, not $0$, we have 
\[K_{\cY_t}+4i(P)(K_{\cG}+\T)\vert_{\cY_t} = g_t^*\big(K_{\cX_t}+4i(P)(K_{\cF_t}+\D_t)\big) + \sum_j b_jF_j.\]
By Lemma~\ref{sep_lc}, we have, for all $j$, that $b_j\geq -1$. 

Moreover, let $X$ be an irreducible component of $\cX_0$ and $Y$ be the proper transform of $X$ in $\cY$. 
Note that $(K_\cG+\Theta)\vert_{Y} - g_0^*(K_\cF+\Delta)\vert_{X}$ is effective since $(\cF,\Delta)$ is a family of stable $\cS_P$-smoothable foliated triples where $g_0$ is the morphism from $Y$ to $X$. 
Moreover, by adjunction, we have 
\begin{align*}
& \Big(\sum_i(a_i-1)E_i + \sum_j b_jF_j\Big)\vert_{Y} \\
= {} & \Big(K_\cY+Y+4i(P)(K_\cG+\T) - g^*\big(K_\cX+X+4i(P)(K_\cF+\D)\big)\Big\vert_{Y} \\
\geq {} & K_{Y} - g_0^*K_{X} \\
\geq {} & -\sum_iE_i\vert_{Y} - \sum_j F_j\vert_{Y}
\end{align*}
where the last inequality follows from $\cX_0$ has semi-log canonical singularities. 
Hence, $a_i\geq 0$ for all $i$. 
\qed
\end{pf}

Now we are ready to prove Theorem~\ref{separatedness}. 

\begin{pf}[Theorem~\ref{separatedness}]
Let $\cY$ be a common resolution of $(\cX_1,\D_1)$, $(\cX_2,\D_2)$, and $\varphi^\circ$. 
Let $g_1$ and $g_2$ be the morphisms fitting into the following commutative diagram. 
\[\xymatrix{ & \cY \ar[ld]_-{g_1} \ar[rd]^-{g_2} & \\
\cX_1 \ar[rd]_-{f_1} \ar@{-->}[rr]^-{\varphi^\circ} & & \cX_2 \ar[ld]^-{f_2} \\
 & T & }\]

Note that $g_1^*\cF_1$ and $g_2^*\cF_2$ are isomorphic over $T^\circ$, so they are isomorphic at the generic point of $\cY$. 
By \cite[Lemma 1.8]{hacon2021birational}, they are in fact isomorphic. 
Let $\cG = g_1^*\cF_1 = g_2^*\cF_2$. 
Moreover, we have the proper transforms $(g_{1*})^{-1}\D_1$ and $(g_{2*})^{-1}\D_2$ are the same. 
Then we put $\T = (g_{1*})^{-1}\D_1 = (g_{2*})^{-1}\D_2$. 

As in subsection~\ref{ss_reduction}, we may assume that $(\cY,\cG,\T) \rw T$ is a semi-stable reduction. 
Note that $\cX_1$ and $\cX_2$ may not be normal after base change. 

Suppose first that $\cX_1$ is normal. 
We write 
\[K_{\cY} + 4i(P)(K_{\cG}+\T) = g_1^*\big(K_{\cX_1}+4i(P)(K_{\cF_1}+\D_1)\big) + \sum_i a_iE_i + \sum_j b_jF_j\]
where $E_i$ are exceptional divisors mapping to $0\in T$ and $F_j$ are flat over $T$. 

Now, by Lemma~\ref{adjunction}, $a_i\geq 0$, and $b_j\geq -1$. 
Also, by adjunction formula, we notice that 
\begin{align*}
& (g_1)_*\left(K_{\cY} + 4i(P)(K_{\cG}+\T)+\sum_j F_j\right) \\ 
= {} & (g_1)_*\left(g_1^*\big(K_{\cX_1}+4i(P)(K_{\cF_1}+\D_1)\big) + \sum_i a_iE_i + \sum_j (b_j+1)F_j\right) \\
= {} & K_{\cX_1}+4i(P)(K_{\cF_1}+\D_1).
\end{align*}
Hence, we have 
\begin{align*}
\cX_1 &\cong \proj\bigoplus_{m\geq 0} H^0(\cX_1,m(K_{\cX_1}+4i(P)(K_{\cF_1}+\D_1))) \\
&= \proj\bigoplus_{m\geq 0} H^0\left(\cY,m\left(K_{\cY} + 4i(P)(K_{\cG}+\T)+\sum_j F_j\right)\right). 
\end{align*}

If $\cX_1$ is not normal, then using the same argument as above, we have the normalization 
\[\cX_1^\nu \cong \proj\bigoplus_{m\geq 0} H^0\left(\cY,m\left(K_{\cY} + 4i(P)(K_{\cG}+\T)+\sum_j F_j\right)\right). \]
Note that the conductor of the normalization $\cX_1^\nu \rw \cX_1$ is in the support of $\sum F_j$. 
Thus we have 
\[\cX_1^\nu \cong \proj\bigoplus_{m\geq 0} H^0\left(\cY,m\left(K_{\cY} + 4i(P)(K_{\cG}+\T)+\sum_j F_j\right)\right) \cong \cX_2^\nu.\]
Therefore $\cX_1$ and $\cX_2$ are isomorphic in codimension 1 where $\cX_1$ and $\cX_2$ here are the original one (before taking a base change).  
Since they both satisfy Serre's condition $S_2$, they are, in fact, isomorphic. 

Now we put $\cX_1\cong\cX_2\cong\cX$.
Notice that $\cF_1$ and $\cF_2$ are isomorphic at the generic point of $\cX$. 
So they are indeed isomorphic. 
Moreover, $\D_1$ and $\D_2$ are isomorphic as well. 
\qed
\end{pf}

\section{Local-closedness for families of foliated triples}
In this section, we will show in Theorem~\ref{open} the local-closedness for families of stable $\cS_P$-smoothable foliated triples. 
(See Definition~\ref{defn_moduli}.)

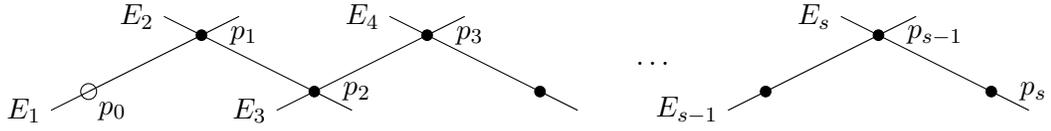
\begin{figure}[h]
\begin{tikzpicture}
\def\r{.5}
\draw (-\r,-0.5*\r) node[left] {$E_1$} -- (0,0) circle (3pt) node[below right] {$p_0$};
\filldraw (0,0) -- (3*\r,1.5*\r) circle (2pt) node[right,xshift=0.25cm] {$p_1$} -- (4*\r,2*\r);
\filldraw (2*\r,2*\r) node[left] {$E_2$} -- (3*\r,1.5*\r) -- (6*\r,0) circle (2pt) node[right,xshift=0.25cm] {$p_2$} -- (7*\r,-0.5*\r);
\filldraw (5*\r,-0.5*\r) node[left] {$E_3$} -- (6*\r,0) -- (9*\r,1.5*\r) circle (2pt) node[right,xshift=0.25cm] {$p_3$} -- (10*\r,2*\r);
\filldraw (8*\r,2*\r) node[left] {$E_4$} -- (9*\r,1.5*\r) -- (12*\r,0*\r) circle (2pt) -- (13*\r,-0.5*\r);
\draw (15*\r,0.75*\r) node {$\ldots$};
\filldraw (17*\r,-0.5*\r) node[left] {$E_{s-1}$} -- (18*\r,0*\r) circle (2pt)  -- (21*\r,1.5*\r) circle (2pt) node[right,xshift=0.25cm] {$p_{s-1}$} -- (22*\r,2*\r);
\filldraw (20*\r,2*\r) node[left] {$E_s$} -- (21*\r,1.5*\r) -- (24*\r,0*\r) circle (2pt) node[right,xshift=0.25cm] {$p_s$} -- (25*\r,-0.5*\r);
\end{tikzpicture}
\caption{A chain of $(-2)$-$\sF$-curves. The straight lines indicate invariant divisors. The circle indicates either a smooth foliation point or a reduced singularity while the solid circles indicate reduced singularities on these curves.}\label{fig_neg2}
\end{figure}

\begin{lem}\label{-2string}
Let $\sF$ be a foliation with only log canonical singularities on a smooth surface $X$. 
Given a Hirzebruch-Jung string $\bigcup_{i=1}^{s}E_i$ of $(-2)$-$\sF$-curves. 
Assume that there are at most two foliation singularities on each $E_i$. 
Then, after re-indexing, we have that $E_1$ has one or two singularities and $E_i$ has exactly two singularities for $i=2$, $\ldots$, $s$. 
Moreover, let $p_{i-1}$ and $p_{i}$ be the singularities on $E_i$ for $i=2$, $\ldots$, $s$ where $p_i$ is the intersection of $E_i$ and $E_{i+1}$ for $i=2$, $\ldots$, $s-1$. 
(See Figure~\ref{fig_neg2}.)
Then we have the following two cases: 
\begin{enumerate}
\item If $E_1$ has only one singularity $p_1$, then $p_1$ is a saddle-node and, for $i\geq 2$, $p_i$ is a reduced and non-degenerate singularity with eigenvalue uniquely determined by $E_2^2$, $\ldots$, $E_s^2$. 
\item If $E_1$ has two singularities $p_0$ and $p_1$, then there is at most one $i\in\{0,\ldots, s\}$ such that $p_i$ is not reduced. 
Moreover, if $\textnormal{Z}(\sF,E_i,p_i)$ (resp. $\textnormal{Z}(\sF,E_{i+1},p_i)$) is non-positive, then $p_{i-1}$ (resp. $p_{i+1}$) is a saddle-node. 
\end{enumerate}
\end{lem}

\begin{pf} 
Since $\sF$ is log canonical, by \cite[Fact I.ii.4]{mcquillan2013almost} the vector field around any fixed point either is a smooth point or belongs to one of three cases in Lemma~\ref{lem_ZCS}. 
\begin{enumerate}
\item If $E_1$ has only one singularity $p_1$, we have that 
\[\textnormal{Z}(\sF,E_1,p_1)=\textnormal{Z}(\sF,E_1)=2.\] 
Thus by Lemma~\ref{lem_ZCS}, we know that $p_1$ is a saddle-node, $\textnormal{Z}(\sF,E_2,p_1) = 1$, and $\textnormal{CS}(\sF,E_2,p_1)=0$. 
By Camacho-Sad formula (Theorem~\ref{thm_CS_formula}), we have that 
\[\lambda_2 := \textnormal{CS}(\sF,E_2,p_2) = E_2^2-\textnormal{CS}(\sF,E_2,p_1) =E_2^2\leq -2. \]
Hence $p_2$ is reduced and non-degenerate. 

\begin{claim}
For all $i=2$, $\ldots$, $s$, $\lambda_i := \textnormal{CS}(\sF,E_{i},p_{i})$ are negative rational numbers less than $-1$ and $p_{i}$ are reduced non-degenerate singularities. 
\end{claim}

We have seen that the claim holds for $i=2$. 
Now we assume that $p_i$'s are reduced and non-degenerate with $\lambda_i\in\bQ_{<-1}$ for $i=2$, $\ldots$, $k-1$. 
Then we have that $\textnormal{Z}(\sF,E_k,p_{k-1})=1$ and $\textnormal{Z}(\sF,E_k,p_k)=1$. 
By Camacho-Sad formula (Theorem~\ref{thm_CS_formula}) and Lemma~\ref{lem_ZCS}, we have that 
\begin{align*}
\lambda_k &= E_k^2-\textnormal{CS}(\sF,E_k,p_{k-1}) \\
&= E_k^2-\frac{1}{\lambda_{k-1}} <  -2+1=-1.
\end{align*}
Therefore, by induction, we complete the proof of the claim and of the first case. 

\item Now we assume that all $E_i$'s have two foliation singularities. 
Suppose there are two non-reduced singularities $p_i$ and $p_j$. 
As the proof of Lemma~\ref{pos_rat}, after a sequence of blowups, we get a string with two separatrices at each singularity, which is impossible by Theorem~\ref{separatrix}. 

Moreover, if $\textnormal{Z}(\sF,E_i,p_i)\leq 0$, then 
\[\textnormal{Z}(\sF,E_i,p_{i-1}) = 2-\textnormal{Z}(\sF,E_i,p_i)\geq 2.\]
So, by Lemma~\ref{lem_ZCS}, we have that $p_{i-1}$ is a saddle-node. 
\qed
\end{enumerate}
\end{pf}

\begin{lem}\label{egl}
Let $\sF$ be a foliation with only log canonical singularities on a smooth surface $X$. 
Given a cycle $E = \bigcup_{i=1}^{s}E_i$ of $(-2)$-$\sF$-curves where $s\geq 2$. 
Then all singularities on $E$ are reduced and non-degenerate. 
\end{lem}
\begin{pf}
\begin{enumerate}
\item If there is only one non-reduced singularity $p$, then we can view this cycle as one Hirzebruch-Jung string of $(-2)$-$\sF$-curves with 2 non-reduced singularities. 
Thus by Lemma~\ref{-2string}, it is impossible. 

\item If there are more than two non-reduced singularities, then we can take any Hirzebruch-Jung string of $(-2)$-$\sF$-curves such that both endpoints are non-reduced. 
Then by Lemma~\ref{-2string}, it is impossible. 

\item Now we know that all singularities on $E$ are reduced. 
Since all $E_i$'s are $(-2)$-$\sF$-curves, we have that all singularities on $E$ are non-degenerate by Lemma~\ref{lem_ZCS}. 
\end{enumerate}
\qed
\end{pf}

\begin{lem}\label{eigenvalue}
Let $\pi:(\cX,\cF) \rw T$ be a smooth family such that 
\begin{enumerate}
\item each irreducible component of foliation singularities is finite over $T$ and pairwise disjoint, and 
\item $\cF_t$ has at worst log canonical singularities for all $t\in T$. 
\end{enumerate}
Given any irreducible component $W$ of singularities of $\cF$, we can define a holomorphic function $\lambda : W \rw \bP^1$ such that $\lambda(w)$ is the eigenvalue of $\cF_{\pi(w)}$ at $w$. 
\end{lem}
\begin{pf}
For any $w\in W$, we can pick an \'{e}tale open subset $\cU$ of $\cX$ containing $w$ such that $(x,y,t)$ are local coordinates of $\cU$ with $\pi(x,y,t) = t$ and $W\cap \cU = (x=y=0)$. 

Then, on $\cU$, the foliation is given by 
\[\big(a(t)x+b(t)y+\textnormal{h.o.t.}\big)\p{x}+\big(c(t)x+d(t)y+\textnormal{h.o.t.}\big)\p{y}\]
where $a$, $b$, $c$, and $d$ are holomorphic functions from $\pi(V)$ to $\bC$ and h.o.t stands for higher order terms (in $x$, $y$, and $t$). 
Thus we form a matrix 
\[A(t) := \left(\begin{matrix}
a(t) & b(t) \\ c(t) & d(t)
\end{matrix}\right)\]
which has two eigenvalues $\mu_+(t)$ and $\mu_-(t)$. 
Since all foliation singularities are log canonical, at least one of two eigenvalues is not zero. 
Therefore, we can define a holomorphic map $\lambda : W\cap\cU \rw \bP^1$ by 
$\lambda(w) = [\mu_+(\pi(w)):\mu_-(\pi(w))]$. 
Notice that this map is independent of the open subset $\cU$ and of change of coordinates. 
Thus, $\lambda$ can be defined on $W$ and is holomorphic. \qed
\end{pf}

\begin{thm}\label{open}
Given a family of foliated triples $f : (\cX,\cF,\D)\rw T$. 
($T$ is not necessarily a curve.) 
Assume that there is a dense subset $T'\subset T$ such that $(\cX,\cF,\D)_t\in\cS_P$ for all $t\in T'$. 
(See Definition~\ref{set_MPdVR} and Definition~\ref{defn_fol_family}.) 
Then there exists an open subset $U\subset T$ such that all fibers of $f$ over $U$ are in $\cS_P$. 
\end{thm}

\begin{pf}
We divide the proof into several steps. 
\begin{enumerate}
\item Let $g: \cY\rw\cX$ be a resolution of singularities, $\cG = g^*\cF$ be the pullback foliation, and $\T$ be the proper transform of $\D$. 
Let $\pi = f\circ g$. 
By Proposition~\ref{ve_Sei}, replacing $\cY$ with its higher model (via some blowups), we may assume that $\cY_t$ is smooth and $\cG_t$ has at worst log canonical singularities for all $t\in U$ for some open subset $U\subset T$. 

\item \textit{Claim.} After shrinking $U$ further if necessary, the following conditions hold: 
\begin{enumerate}
\item $\pi_U$ is smooth. 
\item Each irreducible component of singularities of $\cG$ is finite over $U$ and pairwise disjoint. 
\item Each irreducible component of exceptional divisors of $g$ is flat over $U$. 
\item For any irreducible component $C$ of singularities of $\cG$ and any irreducible component $E$ of exceptional divisors of $g$, we have $C\cap E$ is either $C$ or empty. 
\item For each irreducible component $E$ of exceptional divisors of $g$, either $E_u$ is invariant for all $u\in U$ or $E_u$ is not invariant for all $u\in U$. 
\end{enumerate}
\begin{pf}
More explicitly, (a) comes from the generic smoothness. 
For others, let \[S_1 = \bigcup \pi(C)\] 
where the union is taken over all irreducible components $C$ of singularities of $\cG$ and of exceptional divisors of $g$, which are \emph{not} finite over $U$. 
Let \[S_2 = \bigcup \{t\,\vert\,C_{i,t} = C_{j,t}\}\] 
where the union is taken over all distinct irreducible components $C_i$ and $C_j$ of singularities of $\cG$. 
Let \[S_3 = \bigcup\pi(C\cap E)\] 
where the union is taken over all irreducible components $C$ of singularity of $\cG$ and all irreducible components $E$ of exceptional divisors of $g$ such that $C\cap E$ is neither $C$ nor empty. 

Note that $S_1$, $S_2$, and $S_3$ are all proper closed subsets of $U$. 
So by replacing $U$ by $U\setminus (S_1\cup S_2\cup S_3)$, conditions (a)-(d) hold. 
The condition (e) can be achieved since being non-invariant is an open condition. 
\end{pf}

\item The conditions (a)-(d) will ensure that the number of foliation singularities on an exceptional prime divisor of $g$ is constant over $U$. 
Furthermore, all singularities of $\cG$ do not collide over $U$. 
Also there is a $t_0\in T'\cap U$ since $T'$ is dense in $T$. 

\item We say an irreducible component $E$ of exceptional divisors a \emph{relative $\cG$-exceptional curve} if $E_{t_0}$ is a $\cG_{t_0}$-exceptional curve. 
Then, for any $u\in U$, we have that $E_u$ is invariant, and the contraction of $E_u$ introduces a log canonical foliation singularity on a smooth surface point. 

Contracting all relative $\cG$-exceptional curves, we may assume that there is no relative $\cG$-exceptional curve by induction on Picard number. 

\item\label{G-chain} 
Now we say $E=\bigcup_{i=1}^{s}E_i$ is a \emph{relative $(\cG,\Theta)$-chain} if $E_{t_0} = \bigcup_{i=1}^{s}E_{i,t_0}$ is a $\cG_{t_0}$-chain and $E_{i,t_0}\cdot\T_{t_0}=0$ for $i=1$, $\ldots$, $s$. 
(See Figure~\ref{fig_chain}.)

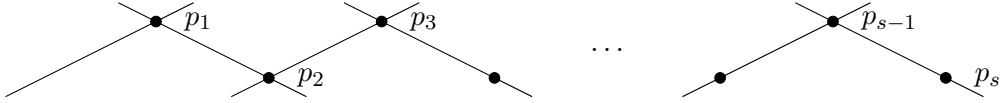
\begin{figure}[h]
\begin{tikzpicture}
\def\r{.5}
\filldraw (-\r,-0.5*\r) -- (0,0) -- (3*\r,1.5*\r) circle (2pt) node[right,xshift=0.25cm] {$p_1$} -- (4*\r,2*\r);
\filldraw (2*\r,2*\r) -- (3*\r,1.5*\r) -- (6*\r,0) circle (2pt) node[right,xshift=0.25cm] {$p_2$} -- (7*\r,-0.5*\r);
\filldraw (5*\r,-0.5*\r) -- (6*\r,0) -- (9*\r,1.5*\r) circle (2pt) node[right,xshift=0.25cm] {$p_3$} -- (10*\r,2*\r);
\filldraw (8*\r,2*\r) -- (9*\r,1.5*\r) -- (12*\r,0*\r) circle (2pt) -- (13*\r,-0.5*\r);
\draw (15*\r,0.75*\r) node {$\ldots$};
\filldraw (17*\r,-0.5*\r) -- (18*\r,0*\r) circle (2pt)  -- (21*\r,1.5*\r) circle (2pt) node[right,xshift=0.25cm] {$p_{s-1}$} -- (22*\r,2*\r);
\filldraw (20*\r,2*\r) -- (21*\r,1.5*\r) -- (24*\r,0*\r) circle (2pt) node[right,xshift=0.25cm] {$p_s$} -- (25*\r,-0.5*\r);
\end{tikzpicture}
\caption{An $\sF$-chain. The straight lines indicate invariant divisors and the solid circles indicate all reduced singularities on these curves.}\label{fig_chain}
\end{figure}

Since, over $t_0$, there is only one singularity $p_{1,t_0}$ on $E_{1,t_0}$, we have that, over any $u\in U$, there is also only one singularity $p_{1,u}$ on $E_{1,u}$. 

By Camacho-Sad formula (Theorem~\ref{thm_CS_formula}), we have $\textnormal{CS}(\cG_u,E_{1,u},p_{1,u}) = E_{1,u}^2\leq -2$. 
And by Lemma~\ref{lem_ZCS}, the eigenvalue $\lambda_1 = \textnormal{CS}(\cG_u,E_{1,u},p_{1,u})$ must be a negative integer less than $-2$. 
Thus, $p_{1,u}$ is a reduced non-degenerate singularity. 

Let the singularities on $E_{k,u}$ be $p_{k-1,u}$ and $p_{k,u}$ for $k = 2$, $\ldots$, $s$. 
\begin{claim}
For all $k=1$, $\ldots$, $s$, $\lambda_k := \textnormal{CS}(\cG_u,E_{k,u},p_{k,u})$ are negative rational numbers less than $-1$ and $p_{k,u}$ are reduced non-degenerate singularities. 
\end{claim}
\begin{pf}
We have shown that this holds when $k=1$. 
Assume the claim holds for $k-1$. 
By Camacho-Sad formula (Theorem~\ref{thm_CS_formula}) and Lemma~\ref{lem_ZCS}, we have that 
\[-2\geq E_{k,u}^2 = \textnormal{CS}(\cG_u,E_{k,u}, p_{k-1,u})+\textnormal{CS}(\cG_u,E_{k,u},p_{k,u}) = \frac{1}{\lambda_{k-1}} + \lambda_{k}. \]
Thus, $\lambda_k\leq -2-\frac{1}{\lambda_{k-1}}<-2+1=-1$. 
Hence, $\lambda_k$ is a negative rational number less than $-1$. 
Therefore, $p_{k,u}$ is also a reduced non-degenerate singularity. 
This completes the proof of the claim. 
\end{pf}

\item We contract all relative $(\cG,\T)$-chains and still using the notation $(\cY,\cG,\T)$ for the resulting foliation. 
Then we have a family such that 
\[(K_\cG+\T)\vert_{\cY_t}\] 
is nef and big for all $t\in T'\cap U$. 
By Lemma~\ref{nef_big_open}, after shrinking $U$ further, we have $(K_\cG+\T)\vert_{\cY_u}$ is nef and big for all $u\in U$. 

Next, in order to show that $(\cY,\cG,\T)_t\in\cS_P$, we will study the \emph{relative $(K_\cG+\T)$-trivial curves}, that is, the irreducible component $E$ of exceptional divisors such that $E_{t_0}$ is $(K_{\cG_{t_0}}+\T_{t_0})$-trivial.  

\item Let $E = E_{0,1}\cup E_{0,2}\cup\bigcup_{i=1}^s E_i$ be the exceptional divisor over a \emph{dihedral singularity} where, for all $u\in U$,  
$\bigcup_{i=1}^sE_{i,u}$ is a Hirzebruch-Jung string, and $E_{0,1,u}$ and $E_{0,2,u}$ are disjoint, of self-intersection $-2$ and both intersect the bad tail $E_{1,u}$. 
(See the Figure~\ref{fig_dual_Dn} for its dual graph.)

\begin{figure}[h]
\begin{tikzpicture}
\def\r{2}
\filldraw (-.707*\r,.707*\r) circle (2pt) node[left] {$E_{0,1}$} -- (0,0) circle (2pt) node[above] {$E_1$}  -- (\r,0) circle (2pt) node[above] {$E_2$};
\filldraw (-.707*\r,-.707*\r) circle (2pt) node[left] {$E_{0,2}$} -- (0,0);
\draw (1.5*\r,0) node {......};
\filldraw (2*\r,0) circle (2pt) node[above] {$E_{s-1}$} -- (3*\r,0) circle (2pt) node[above] {$E_s$};
\end{tikzpicture}
\caption{Dual graph of a $\tn{D}_n$ type singularity}\label{fig_dual_Dn}
\end{figure}
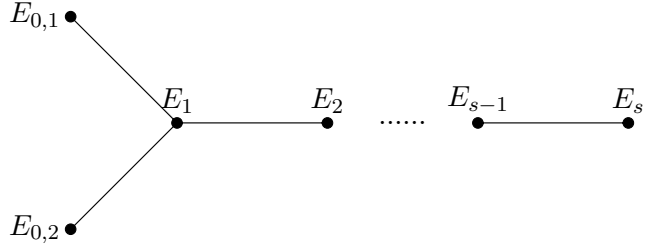

Let $p_{0,1,u}$, $p_{0,2,u}$, and $p_{1,u}$ be three singularities on $E_{1,u}$. 
Let $p_{k-1,u}$ and $p_{k,u}$ be two singularities on $E_{k,u}$ for $k = 2$, $\ldots$, $s$. 
(See the Figure~\ref{fig_Dn}.)

\begin{figure}[h]
\begin{center}
\begin{tikzpicture}
\def\r{.6}
\draw (\r,-0.5*\r) -- (-4*\r,2*\r);
\filldraw (-2*\r,0) -- (-\r,0.5*\r) circle (2pt) node[right,xshift=0.15cm] {$p_{0,2}$} -- (\r,1.5*\r);
\filldraw (-3*\r,0.5*\r) -- (-2*\r,\r) circle (2pt) node[right,xshift=0.15cm] {$p_{0,1}$}  -- (0,2*\r);
\filldraw (-\r,-0.5*\r) -- (0,0) circle (2pt) node[right,xshift=0.15cm] {$p_1$}  -- (3*\r,1.5*\r) circle (2pt) node[right,xshift=0.15cm] {$p_2$} -- (4*\r,2*\r);
\filldraw (2*\r,2*\r) -- (6*\r,0) circle (2pt) node[right,xshift=0.15cm] {$p_3$} -- (7*\r,-0.5*\r);
\filldraw (5*\r,-0.5*\r) -- (9*\r,1.5*\r) circle (2pt) node[right,xshift=0.15cm] {$p_4$} -- (10*\r,2*\r);
\filldraw (8*\r,2*\r) -- (10*\r,\r);
\draw (12*\r,0.75*\r) node {$\ldots\ldots$};
\filldraw (14*\r,\r) -- (15*\r,1.5*\r) circle (2pt) node[right,xshift=0.15cm] {$p_{s-1}$} -- (16*\r,2*\r);
\filldraw (14*\r,2*\r) -- (18*\r,0) circle (2pt) node[right,xshift=0.15cm] {$p_s$} -- (19*\r,-0.5*\r);
\end{tikzpicture}
\caption{ : A graph of $\tn{D}_n$ type. The straight lines indicate invariant divisors and the solid circles indicate all reduced singularities on these divisors.}\label{fig_Dn}
\end{center}
\end{figure}
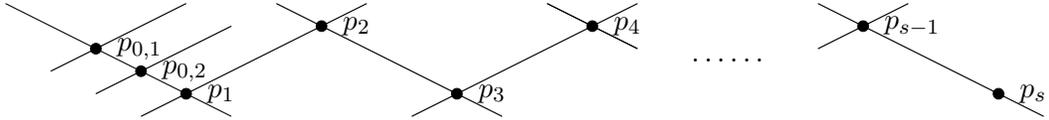

By Camacho-Sad formula (Theorem~\ref{thm_CS_formula}), we have, for any $u\in U$, that 
\[-2 = E_{0,1,u}^2 = \textnormal{CS}(\cG_u,E_{0,1,u},p_{0,1,u}) \mbox{ and } -2 = E_{0,2,u}^2 = \textnormal{CS}(\cG_u,E_{0,2,u},p_{0,2,u}).\]
So $p_{0,1,u}$ and $p_{0,2,u}$ are reduced. 
Thus, 
\begin{align*}
E_{1,u}^2 &= \textnormal{CS}(\cG_u,E_{1,u},p_{1,u}) + \textnormal{CS}(\cG_u,E_{1,u},p_{0,1,u}) + \textnormal{CS}(\cG_u,E_{1,u},p_{0,2,u}) \\
&= \textnormal{CS}(\cG_u,E_{1,u},p_{1,u})-\frac{1}{2}-\frac{1}{2}.
\end{align*}
Hence $\lambda_1 := \textnormal{CS}(\cG_u,E_{1,u},p_{1,u}) = E_{1,u}^2+1\leq -1$.
\begin{claim}
For all $k=1$, $\ldots$, $s$, $\lambda_k := \textnormal{CS}(\cG_u,E_{k,u},p_{k,u})$ are negative rational numbers at most $-1$ and $p_{k,u}$ are reduced non-degenerate singularities. 
\end{claim}
\begin{pf}
We have shown that this holds when $k=1$. 
Assume the claim holds for $k-1$. 
\[-2\geq E_{k,u}^2 = \textnormal{CS}(\cG_u,E_{k,u}, p_{k-1,u})+\textnormal{CS}(\cG_u,E_{k,u},p_{k,u}) = \frac{1}{\lambda_{k-1}} + \lambda_{k}. \]
Thus, $\lambda_k\leq -2-\frac{1}{\lambda_{k-1}}\leq -2+1=-1$. 
Hence, $\lambda_k$ is a negative rational number at most $-1$. 
Therefore, $p_{k,u}$ is also a reduced non-degenerate singularity. 
This completes the proof of the claim. 

So the contraction of $E_u$ introduces a dihedral singularity. 
\end{pf}

\item Let $E = \bigcup_{i=1}^sE_i$ be a connected component of exceptional divisors such that $E_{t_0} = \bigcup_{i=1}^sE_{i,t_0}$ is an \emph{elliptic Gorenstein leaf}. 
Then $E_u = \bigcup_{i=1}^sE_{i,u}$ is also an elliptic Gorenstein leaf for all $u\in U$. 
By Lemma~\ref{egl}, all singularities on $E_u$ are reduced and non-degenerate. 
Therefore, the contraction of $E_u$ introduces a cusp singularity. 

\item (Cyclic non-terminal canonical foliation singularities)
Let $E = \bigcup_{i=1}^sE_i$ be a connected component of exceptional divisors such that the contraction of $E_{t_0} = \bigcup_{i=1}^sE_{i,t_0}$ gives a \emph{cyclic non-terminal} foliation singularity such that $E_{i,t_0}$ are \emph{invariant} for all $i$. 

Then $E_{t_0}$ is a chain of $(-2)$-$\cG_{t_0}$-curves. 
Thus $E_u$ is also a chain of $(-2)$-$\cG_u$-curves for $u\in U$. 
By Lemma~\ref{-2string}, we have the following two cases. 
\begin{enumerate}
\item In the first case, all singularities are reduced. 
So the contraction of $E_u$ introduces a canonical foliation singularity. 
\item In the second case, if all $p_{i,u}$ are reduced (see Figure~\ref{fig_An}), then the contraction of $E_u$ introduces a canonical foliation singularity as well. 

\begin{figure}[h]
\begin{tikzpicture}
\def\r{.5}
\filldraw (-\r,-0.5*\r) -- (0,0) circle (2pt) node[right,xshift=0.25cm] {$p_0$}  -- (3*\r,1.5*\r) circle (2pt) node[right,xshift=0.25cm] {$p_1$} -- (4*\r,2*\r);
\filldraw (2*\r,2*\r) -- (3*\r,1.5*\r) -- (6*\r,0) circle (2pt) node[right,xshift=0.25cm] {$p_2$} -- (7*\r,-0.5*\r);
\filldraw (5*\r,-0.5*\r) -- (6*\r,0) -- (9*\r,1.5*\r) circle (2pt) node[right,xshift=0.25cm] {$p_3$} -- (10*\r,2*\r);
\filldraw (8*\r,2*\r) -- (9*\r,1.5*\r) -- (12*\r,0*\r) circle (2pt) -- (13*\r,-0.5*\r);
\draw (15*\r,0.75*\r) node {......};
\filldraw (17*\r,-0.5*\r) -- (18*\r,0*\r) circle (2pt)  -- (21*\r,1.5*\r) circle (2pt) node[right,xshift=0.25cm] {$p_{s-1}$} -- (22*\r,2*\r);
\filldraw (20*\r,2*\r) -- (21*\r,1.5*\r) -- (24*\r,0*\r) circle (2pt) node[right,xshift=0.25cm] {$p_s$} -- (25*\r,-0.5*\r);
\end{tikzpicture}
\caption{A graph of $\tn{A}_n$ type singularity. The straight lines indicate invariant divisors and the solid circles indicate all reduced singularities on these curves.}\label{fig_An}
\end{figure}
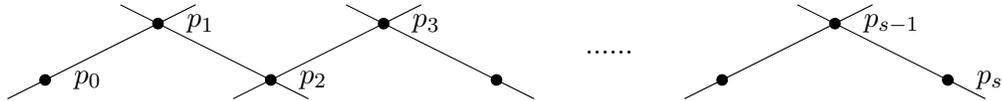

\begin{claim}
If there is an $i$ such that $p_{i,u}$ is not reduced, then all $p_{j,u}$ are non-degenerate for $j\neq i$. 
\end{claim}
\begin{pf}
Suppose not, then one of $p_{i-1,u}$ and $p_{i+1,u}$ is a saddle-node. 
Without loss of generality, we assume that $p_{i+1,u}$ is a saddle-node and 
\[\textnormal{Z}(\cG_u,E_{i+1,u},p_{i,u})\leq 0.\]

Let $W_i$ be the irreducible component of the singularities of $\cG$ containing $p_{i,u}$. 
By Lemma~\ref{eigenvalue}, we have a holomorphic map $\lambda_i : W_i \rw \bP^1$. 
Note that $W_i$ is connected and locally path-connected (in the analytic topology). We have that $W_i$ is path-connected. 
Note that the claim holds true when $\lambda$ is a constant function. 
So we now assume that $\lambda$ is not a constant function. 

There is a sequence $\{w_j\}_{j=1}^\infty$ converging to $p_{i,u}$ such that $\lambda_i(w_j)$ is not a positive rational and converges to $\lambda(p_{i,u})$. 
For any $w_{j}$ close to $p_{i,u}$, there are only two separatrices, one of which is $E_{i+1,\pi(w_j)}$. 
So we may assume $E_{i+1,\pi(w_j)}$ converges to $E_{i+1,u}$. 
Thus $\textnormal{Z}(\cG_u, E_{i+1,u},p_{i,u}) = 1$, which gives a contradiction. 

This completes the proof of the claim. 
\end{pf}

Therefore, the contraction of $E_u$ introduces a mild log canonical foliation singularity. 
\end{enumerate}

\item Let $E = \bigcup_{i=1}^sE_i$ be a connected component of exceptional divisors such that the contraction of $E_{t_0} = \bigcup_{i=1}^sE_{i,t_0}$ gives a \emph{cyclic non-canonical foliation singularity}. 
So there is a unique $j\in\{1, \ldots, s\}$ such that $E_{j,t_0}$ is \emph{non-invariant}. 
Since being non-invariant is an open condition, by shrinking $U$ further, we may assume that $E_{u,t_0}$ is non-invariant for all $u\in U$. 

By~(\ref{G-chain}), when $u$ varies, we know that $\bigcup_{i=1}^{j-1}E_{i,u}$ and $\bigcup_{i=j+1}^{s}E_{i,u}$ have the same structures, respectively. 
Therefore, the contraction of $E_u$ introduces a \emph{mild log canonical} foliation singularities. 

\item Fix any $u\in U$, we may contract all $E_u$  where $E$ is any \emph{relative $(K_\cG+\T)$-trivial curve}. 
It will end up with a \emph{foliated surface} $(Z,\sH)$ with 
\[K_{\cG_u}+\T_u = \pi^*K_\sH\] 
where $\pi : Y \rw Z$ is the contraction morphism. 

\begin{claim}
There is no $K_{\sH}$-trivial curve $C$ with $C^2<0$. 
\end{claim}
\begin{pf}
Suppose $C$ is a $K_\sH$-trivial curve with $C^2<0$. 
Then the proper transform $\wt{C}$ on $Y$ is a $(K_{\cG_u}+\T_u)$-trivial curve. 
Note that $\wt{C}^2\leq C^2<0$. 
Also, by Corollary~\ref{trivial_move}, this $\wt{C}$ moves in a family $\mathscr{C}$ with $\mathscr{C}_u = \wt{C}$. 
However, $\mathscr{C}_{t_0}^2 = \mathscr{C}_{u}^2 = \wt{C}^2<0$ and 
\[(K_{\cG_{t_0}}+\T_{t_0})\cdot\mathscr{C}_{t_0} = (K_{\cG_u}+\T_u)\cdot\mathscr{C}_u = K_\sH\cdot C = 0.\]
Therefore, \emph{$\mathscr{C}$ is a relative $(K_\cG+\T)$-trivial curve}. 
Hence, $\wt{C} = \mathscr{C}_u$ is contracted by $\pi$, which gives a contradiction. 
This completes the proof of the claim. 
\end{pf}

Therefore, $(Z,\sH)$ is a log canonical model and hence $(\cX,\cF,\Delta)_u\in\cS_P$ for any $u\in U$. 
\qed
\end{enumerate}
\end{pf}

\bibliographystyle{amsalpha}
\addcontentsline{toc}{chapter}{\bibname}
\normalem
\bibliography{LogCanMod}

\providecommand{\bysame}{\leavevmode\hbox to3em{\hrulefill}\thinspace}
\providecommand{\MR}{\relax\ifhmode\unskip\space\fi MR }
\providecommand{\MRhref}[2]{%
  \href{http://www.ams.org/mathscinet-getitem?mr=#1}{#2}
}
\providecommand{\href}[2]{#2}
\begin{thebibliography}{KKMSD73}

\bibitem[Bri68]{brieskorn1968rationale}
E.~Brieskorn, \emph{Rationale {S}ingularit\"{a}ten komplexer {F}l\"{a}chen},
  Invent. Math. \textbf{4} (1967/68), 336--358. \MR{222084}

\bibitem[Bru97]{brunella1997feuilletages}
M.~Brunella, \emph{Feuilletages holomorphes sur les surfaces complexes
  compactes}, Ann. Sci. \'{E}cole Norm. Sup. (4) \textbf{30} (1997), no.~5,
  569--594. \MR{1474805}

\bibitem[Bru01]{brunella2001invariance}
\bysame, \emph{Invariance par d\'{e}formations de la dimension de {K}odaira
  d'un feuilletage sur une surface}, Essays on geometry and related topics,
  {V}ol. 1, 2, Monogr. Enseign. Math., vol.~38, Enseignement Math., Geneva,
  2001, pp.~113--132. \MR{1929324}

\bibitem[Bru15]{brunella2015birational}
\bysame, \emph{Birational geometry of foliations}, IMPA Monographs, vol.~1,
  Springer, Cham, 2015. \MR{3328860}

\bibitem[Cam88]{camacho1988quadratic}
C.~Camacho, \emph{Quadratic forms and holomorphic foliations on singular
  surfaces}, Math. Ann. \textbf{282} (1988), no.~2, 177--184. \MR{963011}

\bibitem[CF18]{cascini2018invariance}
P.~Cascini and E.~Floris, \emph{On invariance of plurigenera for foliations on
  surfaces}, J. Reine Angew. Math. \textbf{744} (2018), 201--236. \MR{3871444}

\bibitem[Che21]{chen2021boundedness}
Y.-A. Chen, \emph{Boundedness of minimal partial du {V}al resolutions of
  canonical surface foliations}, Math. Ann. \textbf{381} (2021), no.~1-2,
  557--573. \MR{4322620}

\bibitem[CR14]{cano2014vector}
F.~Cano and C.~Roche, \emph{Vector fields tangent to foliations and blow-ups},
  J. Singul. \textbf{9} (2014), 43--49. \MR{3249046}

\bibitem[CS82]{camacho1982invariant}
C.~Camacho and P.~Sad, \emph{Invariant varieties through singularities of
  holomorphic vector fields}, Ann. of Math. (2) \textbf{115} (1982), no.~3,
  579--595. \MR{657239}

\bibitem[CS20]{cascini2020mmp}
P.~Cascini and C.~Spicer, \emph{On the {MMP} for rank one foliations on
  threefolds}, arXiv preprint arXiv:2012.11433 (2020).

\bibitem[CS21]{cascini2021mmp}
\bysame, \emph{M{MP} for co-rank one foliations on threefolds}, Invent. Math.
  \textbf{225} (2021), no.~2, 603--690. \MR{4285142}

\bibitem[EH87]{eisenbud1987varieties}
D.~Eisenbud and J.~Harris, \emph{On varieties of minimal degree (a centennial
  account)}, Algebraic geometry, {B}owdoin, 1985 ({B}runswick, {M}aine, 1985),
  Proc. Sympos. Pure Math., vol.~46, Amer. Math. Soc., Providence, RI, 1987,
  pp.~3--13. \MR{927946}

\bibitem[Fuj12]{fujino2012minimal}
O.~Fujino, \emph{Minimal model theory for log surfaces}, Publ. Res. Inst. Math.
  Sci. \textbf{48} (2012), no.~2, 339--371. \MR{2928144}

\bibitem[Gro66]{grothendieck1966elements}
A.~Grothendieck, \emph{\'{E}l\'{e}ments de g\'{e}om\'{e}trie alg\'{e}brique.
  {IV}. \'{E}tude locale des sch\'{e}mas et des morphismes de sch\'{e}mas.
  {III}}, Inst. Hautes \'{E}tudes Sci. Publ. Math. (1966), no.~28, 255.
  \MR{217086}

\bibitem[Har77]{hartshorne1977algebraic}
R.~Hartshorne, \emph{Algebraic geometry}, Springer-Verlag, New York-Heidelberg,
  1977, Graduate Texts in Mathematics, No. 52. \MR{0463157}

\bibitem[HK10]{hacon2010classification}
C.~D. Hacon and S.~J. Kov\'{a}cs, \emph{Classification of higher dimensional
  algebraic varieties}, Oberwolfach Seminars, vol.~41, Birkh\"{a}user Verlag,
  Basel, 2010. \MR{2675555}

\bibitem[HL21]{hacon2021birational}
C.~D. Hacon and A.~Langer, \emph{On birational boundedness of foliated
  surfaces}, J. Reine Angew. Math. \textbf{770} (2021), 205--229. \MR{4193468}

\bibitem[KKMSD73]{kempf1973toroidal}
G.~Kempf, F.~Knudsen, D.~Mumford, and B.~Saint-Donat, \emph{Toroidal
  embeddings. {I}}, Lecture Notes in Mathematics, Vol. 339, Springer-Verlag,
  Berlin-New York, 1973. \MR{0335518}

\bibitem[KM98]{kollar1998birational}
J.~Koll\'{a}r and S.~Mori, \emph{Birational geometry of algebraic varieties},
  Cambridge Tracts in Mathematics, vol. 134, Cambridge University Press,
  Cambridge, 1998, With the collaboration of C. H. Clemens and A. Corti,
  Translated from the 1998 Japanese original. \MR{1658959}

\bibitem[Kol08]{kollar2008hulls}
J.~Koll\'{a}r, \emph{Hulls and husks}, arXiv e-print, arXiv:0805.0576 (2008).

\bibitem[Lan00]{langer2000chern}
A.~Langer, \emph{Chern classes of reflexive sheaves on normal surfaces}, Math.
  Z. \textbf{235} (2000), no.~3, 591--614. \MR{1800214}

\bibitem[Lan01]{langer2001adjoint}
\bysame, \emph{Adjoint linear systems on normal log surfaces}, Compositio Math.
  \textbf{129} (2001), no.~1, 47--66. \MR{1856022}

\bibitem[Laz04]{lazarsfeld2004positivity}
R.~Lazarsfeld, \emph{Positivity in algebraic geometry. {I}}, Ergebnisse der
  Mathematik und ihrer Grenzgebiete. 3. Folge. A Series of Modern Surveys in
  Mathematics [Results in Mathematics and Related Areas. 3rd Series. A Series
  of Modern Surveys in Mathematics], vol.~48, Springer-Verlag, Berlin, 2004,
  Classical setting: line bundles and linear series. \MR{2095471}

\bibitem[LN88]{linsneto1988algebraic}
A.~Lins~Neto, \emph{Algebraic solutions of polynomial differential equations
  and foliations in dimension two}, Holomorphic dynamics ({M}exico, 1986),
  Lecture Notes in Math., vol. 1345, Springer, Berlin, 1988, pp.~192--232.
  \MR{980960}

\bibitem[McQ08]{mcquillan2008canonical}
M.~McQuillan, \emph{Canonical models of foliations}, Pure Appl. Math. Q.
  \textbf{4} (2008), no.~3, Special Issue: In honor of Fedor Bogomolov. Part 2,
  877--1012. \MR{2435846}

\bibitem[MM80]{mattei1980holonomie}
J.-F. Mattei and R.~Moussu, \emph{Holonomie et int\'{e}grales premi\`eres},
  Ann. Sci. \'{E}cole Norm. Sup. (4) \textbf{13} (1980), no.~4, 469--523.
  \MR{608290}

\bibitem[MP13]{mcquillan2013almost}
M.~McQuillan and D.~Panazzolo, \emph{Almost \'{e}tale resolution of
  foliations}, J. Differential Geom. \textbf{95} (2013), no.~2, 279--319.
  \MR{3128985}

\bibitem[PS19]{pereira2019effective}
J.~V. Pereira and R.~Svaldi, \emph{Effective algebraic integration in bounded
  genus}, Algebr. Geom. \textbf{6} (2019), no.~4, 454--485. \MR{3957403}

\bibitem[Rei]{reid2012surface}
M.~Reid, \emph{Surface cyclic quotient singularities and {H}irzebruch-{J}ung
  resolutions}, manuscript available at
  \url{http://homepages.warwick.ac.uk/~masda/surf/more/cyclic.pdf}, Accessed:
  2022-02-24.

\bibitem[Rei87]{reid1987canonical}
\bysame, \emph{Young person's guide to canonical singularities}, Algebraic
  geometry, {B}owdoin, 1985 ({B}runswick, {M}aine, 1985), Proc. Sympos. Pure
  Math., vol.~46, Amer. Math. Soc., Providence, RI, 1987, pp.~345--414.
  \MR{927963}

\bibitem[Sak84]{sakai1984weil}
F.~Sakai, \emph{Weil divisors on normal surfaces}, Duke Math. J. \textbf{51}
  (1984), no.~4, 877--887. \MR{771385}

\bibitem[Seb97]{sebastiani1997existence}
M.~Sebastiani, \emph{Sur l'existence de s\'{e}paratrices locales des
  feuilletages des surfaces}, An. Acad. Brasil. Ci\^{e}nc. \textbf{69} (1997),
  no.~2, 159--162. \MR{1754036}

\bibitem[Sei68]{seidenberg1968reduction}
A.~Seidenberg, \emph{Reduction of singularities of the differential equation
  {$A\,dy=B\,dx$}}, Amer. J. Math. \textbf{90} (1968), 248--269. \MR{220710}

\bibitem[Sim94]{simpson1994moduli}
C.~T. Simpson, \emph{Moduli of representations of the fundamental group of a
  smooth projective variety. {I}}, Inst. Hautes \'{E}tudes Sci. Publ. Math.
  (1994), no.~79, 47--129. \MR{1307297}

\bibitem[Spi20]{spicer2020higher}
C.~Spicer, \emph{Higher-dimensional foliated {M}ori theory}, Compos. Math.
  \textbf{156} (2020), no.~1, 1--38. \MR{4036447}

\bibitem[Suw98]{suwa1998indices}
T.~Suwa, \emph{Indices of vector fields and residues of singular holomorphic
  foliations}, Actualit\'{e}s Math\'{e}matiques. [Current Mathematical Topics],
  Hermann, Paris, 1998. \MR{1649358}

\end{thebibliography}

\end{document}